\documentclass[12pt]{amsart}

\usepackage{graphicx,amssymb,latexsym,amsfonts,txfonts,amsmath,amsthm}
\usepackage{pdfsync,color,tabularx,rotating}
\usepackage[all,cmtip]{xy}
\usepackage{url}
\usepackage{tikz}
\usepackage{amssymb}
\usepackage{comment}

\advance\hoffset by -0.9 truecm

\newtheorem{theo}{Theorem}[section]
\newtheorem{prop}[theo]{Proposition}
\newtheorem{coro}[theo]{Corollary}
\newtheorem{lemm}[theo]{Lemma}

\theoremstyle{remark}

\theoremstyle{remark}

\theoremstyle{remark}

\theoremstyle{remark}

\begin{document}

\title{Regular dessins with primitive automorphism groups}

\author{Gareth A. Jones and Martin Ma\v caj }

\address{School of Mathematical Sciences, University of Southampton, Southampton SO17 1BJ, UK}
\email{G.A.Jones@maths.soton.ac.uk}

\address{Faculty of Mathematics, Physics and Informatics, Comenius University, Bratislava, Slovakia}

\email{macaj@dcs.fmph.uniba.sk
}



\begin{abstract}
We classify the dessins $\mathcal D$ for which the automorphism group $G$ acts primitively and faithfully on the points over one of the three critical values (without loss of generality the black vertices in the usual bipartite map representation). We show that they are all generalised Paley dessins, in which the black vertices are the elements of a finite field ${\mathbb F}_q$, and $G$ is a subgroup of the affine group ${\rm AGL}_1(q)$. Using earlier results obtained with Streit and Wolfart, we determine the orbits of the absolute Galois group on these dessins, we show that they are all defined over certain cyclotomic fields, and we obtain defining equations in some special cases. Relaxing the condition of a faithful action allows only cyclic regular coverings of these dessins.
\end{abstract}

\maketitle


\section{Introduction}\label{sec:intro}

\bigskip

Every dessin is a quotient of a regular dessin (one with the greatest possible symmetry) by a subgroup of its automorphism group, so much of the study of dessins can be reduced to that of regular dessins and their automorphism groups. Unfortunately, the great proliferation of regular dessins (resulting from that of two-generator finite groups) makes it almost impossible to study them systematically. They all have the property that the automorphism group acts transitively on each of the three sets of critical points, the fibres over the critical values of the Bely\u\i\/ function. However, if we strengthen this transitivity property by requiring that the action on at least one of these sets should be primitive (preserving no nontrivial equivalence relations), then only a very specific and easily-understood family of dessins arises, namely the generalised Paley dessins (when the action is also faithful) and their cyclic regular covers (when it is not faithful).

We will say that a dessin $\mathcal D$ is {\sl primitive\/} or {\sl faithful\/} if its automorphism group $G={\rm Aut}\,{\mathcal D}$ acts primitively or faithfully on the set of black vertices. (Throughout this paper we will represent dessins as bipartite graphs on surfaces, with vertices coloured black or white; by applying suitable dualities of dessins one can replace statements about black vertices with similar statements about white vertices or face centres.) 
Our main result is as follows:

\begin{theo}\label{theo:main}
A regular dessin is primitive and faithful if and only if it is a generalised Paley dessin.
\end{theo}

In a generalised Paley dessin $\mathcal D$, the black vertices are identified with the elements of a finite field ${\mathbb F}_q$ of order $q=p^d$ for some prime $p$, and the automorphism group $G={\rm Aut}\,{\mathcal D}$ is a subgroup of the $1$-dimensional affine group ${\rm AGL}_1(q)$ acting primitively on the field. Equivalently, $G$ is a semidirect product $T\rtimes S$, where $T$ is the translation subgroup (an elementary abelian group of order $q$) and $S$ is a subgroup of the multiplicative group ${\mathbb F}_q^*={\mathbb F}_q\setminus\{0\}$ of ${\mathbb F}_q$ acting on the field as an irreducible subgroup of ${\rm GL}_d(p)$, so that $d$ is the multiplicative order of $p$ mod~$(n)$ where $n=|S|$. The standard generators of $G$ are $x$, $y$ and $z$, where $x$ is a generator of $S$, $y$ is an arbitrary element of $G\setminus S$, and $z=(xy)^{-1}$. Thus the black vertices have valency $n$, while the valencies of the white vertices and the faces, and hence the type and genus of the dessin, depend on which coset $Tx^c$ of $T$ contains $y$. In particular, if $q\equiv 1$ mod~$(4)$ with $n=(q-1)/2$ and $c=n/2$, then $S$ is the group of quadratic residues in ${\mathbb F}_q$, and the white vertices have valency $2$, so that if they are ignored the embedded graph is the Paley graph $P_q$, in which vertices are adjacent if and only if they differ by a quadratic residue.

For the origins of the Paley graphs, see~\cite{Jones2}. Generalised Paley graphs were introduced by Lim and Praeger in~\cite{LP}, and generalised Paley maps were studied by Jones in~\cite{Jones}. These maps played an important role in the classification by Jones and Ma\v caj in~\cite{JM} of regular maps with primitive automorphism groups. Generalised Paley dessins can be regarded as their natural generalisation to orientably regular hypermaps.

For a given $n$ and prime $p$ not dividing $n$, the number of non-isomorphic generalised Paley dessins $\mathcal D$ is $n\phi(n)/d$ where $\phi$ is Euler's function, with $\phi(n)/d$ dessins for each choice of $c\in{\mathbb Z}_n$. These $\phi(n)/d$ dessins are permuted transitively by the group of hole operations $H_j$ for $j$ coprime to $n$. It follows from earlier results of Jones, Streit and Wolfart~\cite{JSW} that they form orbits under the absolute Galois group ${\rm Gal}\,\overline{\mathbb Q}/{\mathbb Q}$, and are defined over a specific subfield of the cyclotomic field of the $n$th roots of unity, namely the splitting field of~$p$. This allows us to give defining equations for a few special cases of these dessins.

After two sections giving background on dessins and some simple examples, the basic properties of generalised Paley dessins are explained and determined in Sections~\ref{sec:Paley}--\ref{sec:quotients}, and further examples are given in Section~\ref{sec:further}, with some of the techniques illustrated for $n=7$ and $p=29$ in Section~\ref{sec:calculations}. The main theorem is proved in Section~\ref{sec:proof}. Galois conjugacy and defining equations are considered in Sections~\ref{sec:Galois} and~\ref{sec:equations}, and there is a brief discussion of non-faithful primitive dessins in Section~\ref{sec:non-faithful}. A technical problem concerning number theory is dealt with in Section~\ref{sec:appendix}.


\section{Dessins d'enfants}\label{sec:dessins}

This is just a brief outline of the theory of dessins; for more details see~\cite{JW}.
By Bely\u\i's Theorem, a compact Riemann surface $\mathcal S$ is defined, as an algebraic curve, over an algebraic number field if and only if it admits a Bely\u\i\/ function, that is, a non-constant meromorphic function $\beta:{\mathcal S}\to\Sigma={\mathbb P}^1({\mathbb C})={\mathbb C}\cup\{\infty\}$ with at most three critical values. Applying a M\"obius transformation to $\Sigma$ if necessary, one can assume these are $0$, $1$ and $\infty$. Then $\beta$ can be represented combinatorially by a {\sl dessin d'enfant} $\mathcal D$, a bipartite map on $\mathcal S$ with black and white vertices at the points over $0$ and $1$, and edges over the unit interval $[0,1]\subset{\mathbb R}$; the edges correspond bijectively to the sheets of the covering $\beta$, and the faces correspond to the points over $\infty$, called the face-centres. (It is sometimes useful to extend $\mathcal D$ to a tripartite map, with red vertices over $\infty$ and edges over $\mathbb R$.)

The automorphism group $G={\rm Aut}\,{\mathcal D}$ of $\mathcal D$ is the group of covering transformations of $\beta$, acting semi-regularly on the edges. We say that $\mathcal D$ is {\sl regular\/} if $\beta$ is a regular covering, or equivalently $G$ acts regularly on the edges. In this case $G$ acts transitively on the black and the white vertices and the face-centres, so all black and all white vertices have the same valencies $n$ and $m$, and all faces are $2l$-gons for some $l$. The triple $(n,m,l)$ is called the {\sl type\/} of $\mathcal D$. Then $G$ has generators $x, y, z$ (rotations around an incident black and white vertex and face-centre) satisfying
\[x^n=y^m=z^l=xyz=1.\]

If $m\le 2$ one can ignore the white vertices of $\mathcal D$ (or more precisely, regard them as ordinary points on their incident edge or edges), giving an oriented map ${\mathcal M}={\mathcal M}({\mathcal D})$ in which the vertices are the black vertices of $\mathcal D$, the face are the same, and the edges are the unions of the edges incident with each white vertex of $\mathcal D$. A second way of converting a dessin into a map, applicable to every dessin $\mathcal D$, is to ignore vertex colours, regarding $\mathcal D$ as an uncoloured map ${\mathcal D}^{\circ}$ on the same surface.


\section{Simple examples}\label{sec:examples}

First we consider some simple examples of primitive and faithful actions.

\medskip

\begin{figure}[h!]
\label{fig:cube}
\begin{center}
\begin{tikzpicture}[scale=0.15, inner sep=0.7mm]

\node (a) at (0,0) [shape=circle, draw] {};

\node (b) at (10,0)  [shape=circle, fill=black]  {};
\node (c) at (3,9.5)  [shape=circle, fill=black] {};
\node (d) at (-8,5.9) [shape=circle, fill=black] {};
\node (e) at (-8,-5.9) [shape=circle, fill=black] {};
\node (f) at (3,-9.5) [shape=circle, fill=black] {};

\draw [thick] (a) to (b);
\draw [thick] (a) to (c);
\draw [thick] (a) to (d);
\draw [thick] (a) to (e);
\draw [thick] (a) to (f);

\end{tikzpicture}

\end{center}
\caption{The star dessin ${\mathcal St}_5$}
\label{fig:St5}
\end{figure}
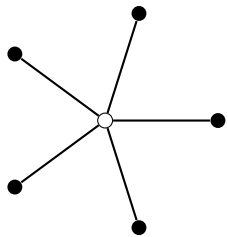

\noindent{\bf Example 1} Let $\mathcal D$ be the regular planar dessin of type $(1,m,m)$ embedding a complete bipartite graph $K_{m,1}$ for some integer $m\ge 1$. This has one white vertex, joined by a single edge to each of the $m$ black vertices. We will call $\mathcal D$ an $m$-{\sl star dessin\/} and denote it by ${\mathcal St}_m$. The corresponding Bely\u\i\/ function $\Sigma\to\Sigma$ is $z\mapsto 1-z^m$. The automorphism group $G$ is a cyclic group ${\rm C}_m$ of order $m$, generated by $z\mapsto e^{2\pi i/m}z$. This group acts transitively and faithfully on the black vertices (the $m$th roots of $1$ in $\mathbb C$), and acts primitively on them if and only if $m$ is prime. Figure~\ref{fig:St5} shows ${\mathcal St}_5$.

\medskip

\begin{lemm}\label{lem:regonV}
Let $\mathcal D$ be a regular dessin with automorphism group $G$. Then $G$ acts regularly (and thus faithfully) on the black vertices of $\mathcal D$ if and only if ${\mathcal D}\cong{\mathcal St}_m$ for some $m\ge 1$.
\end{lemm}

\noindent{\sl Proof.} The star dessins clearly have this property. For the converse, if $G$ acts regularly on the black vertices of $\mathcal D$, then since $G$ also acts regularly on the edges, there must be the same number of edges as black vertices, and hence the latter must all have valency $1$; since the embedded graph is connected there can therefore be only one white vertex, so this graph is a complete bipartite graph $K_{m,1}$ for some $m$. Up to isomorphism, this bicoloured graph has a unique embedding as a dessin, namely as ${\mathcal St}_m$.
\hfill$\square$

\medskip

\begin{figure}[h!]

\begin{center}
 \begin{tikzpicture}[scale=0.6, inner sep=0.8mm]

\node (A) at (-12.5,-1.73) [shape=circle, draw, fill=black] {};
\node (B) at (-10.5,-1.73) [shape=circle, draw] {};
\node (C) at (-9.5,0) [shape=circle, fill=black] {}; 
\node (D) at (-10.5,1.73) [shape=circle, draw] {};
\node (E) at (-12.5,1.73) [shape=circle, fill=black] {};
\node (F) at (-13.5,0) [shape=circle, draw] {};
\draw (A) to (B);
\draw (B) to (C);
\draw (C) to (D);
\draw (D) to (E);
\draw (E) to (F);
\draw (F) to (A);

\draw[dashed] (-14.5,-1.73) to (-11.5,-3.54);
\draw[dashed] (-11.5,-3.54) to (-8.5,-1.73);
\draw[dashed] (-8.5,-1.73) to (-8.5,1.73);
\draw[dashed] (-8.5,1.73) to (-11.5,3.46);
\draw[dashed] (-11.5,3.46) to (-14.5,1.73);
\draw[dashed] (-14.5,1.73) to (-14.5,-1.73);

\draw (A) to (-13,-2.63);
\draw (B) to (-10,-2.63);
\draw (C) to (-8.5,0);
\draw (D) to (-10,2.6);
\draw (E) to (-13,2.6);
\draw (F) to (-14.5,0);

\node at (-16,0) {${\mathcal F}_3$};


\node (A) at (2,0) [shape=circle, draw, fill=black] {};
\node (B) at (1.41,1.41) [shape=circle, draw] {};
\node (C) at (0,2) [shape=circle, fill=black] {}; 
\node (D) at (-1.41,1.41) [shape=circle, draw] {};
\node (E) at (-2,0) [shape=circle, fill=black] {};
\node (F) at (-1.41,-1.41) [shape=circle, draw] {};
\node (G) at (0,-2) [shape=circle, fill=black] {};
\node (H) at (1.41,-1.41) [shape=circle, draw] {};
\draw (A) to (B);
\draw (B) to (C);
\draw (C) to (D);
\draw (D) to (E);
\draw (E) to (F);
\draw (F) to (G);
\draw (G) to (H);
\draw (H) to (A);

\draw[dashed] (3,0) to (2.77,1.14);
\draw[dashed] (2.77,1.14) to (2.12,2.12);
\draw[dashed] (2.12,2.12) to (1.14,2.77);
\draw[dashed] (1.14,2.77) to (0,3);

\draw[dashed] (3,0) to (2.77,-1.14);
\draw[dashed] (2.77,-1.14) to (2.12,-2.12);
\draw[dashed] (2.12,-2.12) to (1.14,-2.77);
\draw[dashed] (1.14,-2.77) to (0,-3);

\draw[dashed] (-3,0) to (-2.77,1.14);
\draw[dashed] (-2.77,1.14) to (-2.12,2.12);
\draw[dashed] (-2.12,2.12) to (-1.14,2.77);
\draw[dashed] (-1.14,2.77) to (0,3);

\draw[dashed] (-3,0) to (-2.77,-1.14);
\draw[dashed] (-2.77,-1.14) to (-2.12,-2.12);
\draw[dashed] (-2.12,-2.12) to (-1.14,-2.77);
\draw[dashed] (-1.14,-2.77) to (0,-3);

\draw (A) to (2.87,0.57);
\draw (A) to (2.87,-0.57);
\draw (C) to (0.57,2.87);
\draw (C) to (-0.57,2.87);
\draw (E) to (-2.87,0.57);
\draw (E) to (-2.87,-0.57);
\draw (G) to (0.57,-2.87);
\draw (G) to (-0.57,-2.87);

\draw (B) to (2.45,1.63);
\draw (B) to (1.63,2.45);
\draw (D) to (-2.45,1.63);
\draw (D) to (-1.63,2.45);
\draw (F) to (-2.45,-1.63);
\draw (F) to (-1.63,-2.45);
\draw (H) to (2.45,-1.63);
\draw (H) to (1.63,-2.45);

\node at (0.7,3.4) {$A$};
\node at (2.8,-1.8) {$A$};

\node at (-4.5,0) {${\mathcal F}_4$};

\end{tikzpicture}

\end{center}
\caption{The Fermat dessins ${\mathcal F}_3$ and ${\mathcal F}_4$.}
\label{fig:F3F4} 
\end{figure}
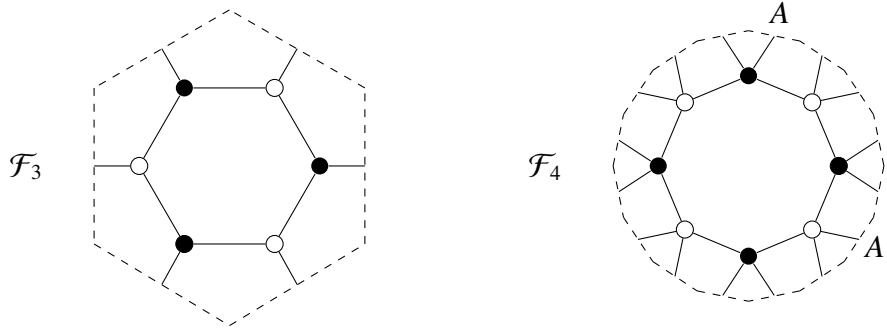

\noindent{\bf Example 2} The Fermat dessin ${\mathcal F}_n$ of degree $n$ is regular, of type $(n,n,n)$ and genus $(n-1)(n-2)/2$, with
\[G=\langle x,y\mid x^n=y^n=[x,y]=1\rangle\cong{\rm C}_n\times{\rm C}_n.\]
The embedded graph is the complete bipartite graph $K_{n,n}$. The dessins ${\mathcal F}_3$ and  ${\mathcal F}_4$ are shown in Figure~\ref{fig:F3F4}: in the case of ${\mathcal F}_3$ opposite sides of the outer hexagon are identified to form a torus, while in the case of ${\mathcal F}_4$ sides of the outer $16$-gon are identified, $A$ with $A$ and the others using the obvious ${\rm D}_4$ symmetry, to give a surface of genus $3$.

The automorphism group $G$
acts transitively but not faithfully on the $n$ black vertices, on the $n$ white vertices and on the $n$ faces. The kernels of these actions of $G$ are the cyclic subgroups generated by the canonical generators $x$, $y$ and $z:=(xy)^{-1}$, so in each case $G$ induces a permutation group $\overline{G}\cong{\rm C}_n$ acting regularly. Each of these three actions is primitive if and only if $n$ is prime. The kernels of these three actions of $G$ have trivial intersection, so the intransitive action of $G$ on the set of all $3n$ critical points is faithful.

\medskip

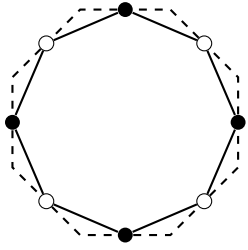
\begin{figure}[h!]
\label{fig:cube}
\begin{center}
\begin{tikzpicture}[scale=0.15, inner sep=0.7mm]

\node (a) at (10,0)  [shape=circle, fill=black]  {};
\node (b) at (7,7) [shape=circle, draw] {} ;
\node (c) at (0,10)  [shape=circle, fill=black] {};
\node (d) at (-7,7 )[shape=circle, draw] {};
\node (e) at (-10,0) [shape=circle, fill=black] {};
\node (f) at (-7,-7) [shape=circle, draw] {};
\node (g) at (0,-10)  [shape=circle, fill=black] {};
\node (h) at (7,-7) [shape=circle, draw] {};

\draw [thick] (a) to (b) to (c) to (d) to (e) to (f) to (g) to (h) to (a);

\draw [thick, dashed] (a) to (10,4) to (b) to (4,10) to (c) to (-4,10) to (d) to (-10,4) to (e) to (-10,-4) to (f) to (-4,-10) to (g) to (4,-10) to (h) to (10,-4) to (a);

\end{tikzpicture}

\end{center}
\caption{The dessin $\mathcal Q$ (identify opposite sides of the outer octagon).}
\label{fig:Q8}
\end{figure}

\noindent{\bf Example 3} The quaternion group ${\rm Q}_8$ of order $8$ has, up to automorphisms, a unique generating triple of elements $x, y, z$ of order $4$ satisfying $xyz=1$. They determine a regular dessin $\mathcal Q$ of type $(4,4,4)$ and genus $2$ with automorphism group $G\cong {\rm Q}_8$. This is shown in Figure~\ref{fig:Q8} where opposite sides of the outer hexagon are identified with each other. There are two black vertices, two white vertices, and two faces, each pair transposed by $G$. In each case the action is primitive, and its kernel is the commutator subgroup $G'$ (also the centre) of $G$, of order $2$, generated by a half turn about the centre of the diagram.
The corresponding hypermap is RPH2.4 in Conder's census of regular hypermaps~\cite{Conder}. The quotient ${\mathcal Q}/G'$ is the Fermat dessin ${\mathcal F}_2$, an embedding of $K_{2,2}$ in the sphere.

\medskip

\noindent{\bf Example 4} The preceding example can be extended to all odd primes, though the structure of the automorphism group is slightly different. For any prime $p>2$, let
\[G=\langle a, b \mid a^{p^2}=b^p=1,\, a^b=a^{p+1}\rangle \cong {\rm C}_{p^2}\rtimes {\rm C}_p\]
be the unique nonabelian group of order $p^3$ and exponent $p^2$. Each element has the unique form $a^ib^j$ for $i\in{\mathbb Z}_{p^2}$ and $j\in{\mathbb Z}_p$, and it has order $p^2$ if and only if $i$ is coprime to $p$. The commutator subgroup $G'$ (also the centre) is $\langle a^p\rangle\cong{\rm C}_p$, with $G/G'\cong {\rm C}_p\times{\rm C}_p$. By using standard generators $x=a$, $y=ab$ and $z=(xy)^{-1}=a^{-2(p+1)}b^{-1}$ of order $p^2$ we obtain a regular dessin ${\mathcal D}_p$ of type $(p^2, p^2, p^2)$ and genus $(p-1)^2(p+2)/2$, with automorphism group $G$ acting primitively on the $p$ black vertices, on the $p$ white vertices, and on the $p$ faces. The kernels of these actions, namely the normal subgroups $\langle x\rangle$, $\langle y\rangle$ and $\langle z\rangle$ of index $p$ in $G$, all contain $G'$, so the action of $G$ on the $3p$ special points of ${\mathcal D}_p$ is not faithful, with kernel $G'$. This dessin ${\mathcal D}_p$ is a regular $p$-sheeted covering of the Fermat dessin ${\mathcal D}_p/A'\cong{\mathcal F}_p$ of degree $p$, branched over all $3p$ of its critical points; for $p=3$ or $5$ the corresponding hypermap is RPH10.32 or RPH56.70 in~\cite{Conder}. 


\section{Generalised Paley dessins}\label{sec:Paley}

For further examples, and an important definition, we need the following lemma:

\begin{lemm}\label{lem:affineprim}
Let $G$ be a subgroup of the affine group ${\rm AGL}_1(q)$ for some prime power $q=p^d$. Then $G$ acts primitively on the field ${\mathbb F}_q$ if and only if 
\[G=G_S:=\{t\mapsto at+b\mid a\in S, \, b\in{\mathbb F}_q\}\]
for some subgroup $S$ of order $n$ in the multiplicative group ${\mathbb F}_q^*={\mathbb F}_q\setminus\{0\}$, where either $n=d=1$, or $n>1$ and $S$ satisfies the following equivalent conditions:
\begin{itemize}
\item [(a)] $S$ acts by multiplication on ${\mathbb F}_q$ as an irreducible subgroup of ${\rm GL}_d(p)$;
\item [(b)] $S$ generates the additive group of ${\mathbb F}_q$;
\item [(c)] $p$ has multiplicative order $d$ {\rm mod}~$(n)$.
\end{itemize}
\end{lemm}
The proof is elementary, the only point requiring comment being that if $G$ is primitive then it is transitive on ${\mathbb F}_q$, so it has order divisible by $q$ and hence contains the unique Sylow $p$-subgroup of ${\rm AGL}_1(q)$, namely the translation subgroup $T=\{t\mapsto t+b\mid b\in{\mathbb F}_q\}$ of order $q$. Note that $G_S$ is a semidirect product $T\rtimes S$ of $T$, acting regularly on ${\mathbb F}_q$ by translations, and the stabiliser $S=G_0$ of $0$, acting on ${\mathbb F}_q$ by multiplication. In the case when $n=d=1$ we have $G_S=T\cong{\rm C}_p$.

\medskip

We can use this result to construct a wider class of primitive dessins. Let $G=G_S$ where $S$ satisfies the conditions of Lemma~\ref{lem:affineprim}. Since ${\mathbb F}_q^*$ is cyclic, so is its subgroup $S$, so let $x$ be any generator of $S$. If $y$ is any element of $G\setminus S$, then $G=\langle x, y\rangle$ since $S$, as a point stabiliser in the primitive group $G$, is a maximal subgroup of $G$. We can therefore define a regular dessin $\mathcal D$ with ${\rm Aut}\,{\mathcal D}\cong G$ by taking the generating triple $(x, y, z)$ for $G$ where $z:=(xy)^{-1}$. The black vertices of $\mathcal D$, of valency $n$, are identified with the cosets of $\langle x\rangle=S$ in $G$, and hence with the elements of ${\mathbb F}_q$, in each case admitting the natural action of $G$, so $G$ acts primitively and faithfully on them. We will call such a dessin $\mathcal D$ a {\sl generalised Paley dessin\/}, since this construction is the natural extension to dessins of the concept of a generalised Paley map (see~\cite{Jones} or \cite[Section 9.1.1]{JW}).

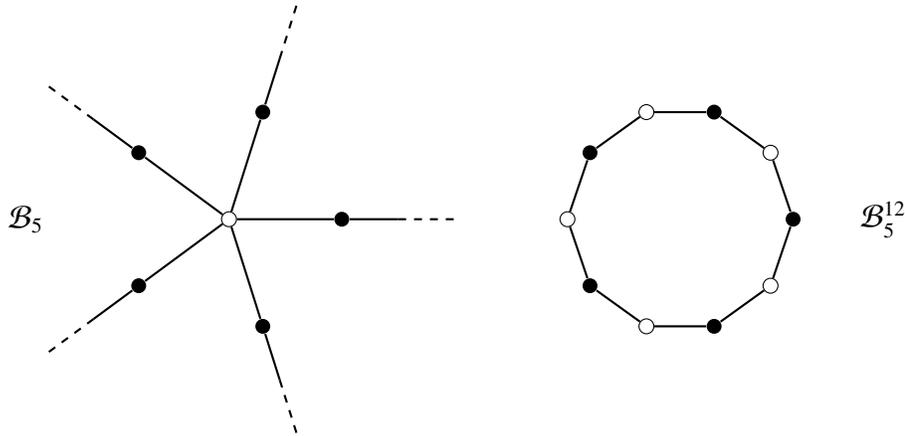
\begin{figure}[h!]
\label{fig:cube}
\begin{center}
\begin{tikzpicture}[scale=0.15, inner sep=0.7mm]

\node (a) at (0,0) [shape=circle, draw] {};

\node (b) at (10,0)  [shape=circle, fill=black]  {};
\node (c) at (3,9.5)  [shape=circle, fill=black] {};
\node (d) at (-8,5.9) [shape=circle, fill=black] {};
\node (e) at (-8,-5.9) [shape=circle, fill=black] {};
\node (f) at (3,-9.5) [shape=circle, fill=black] {};

\draw [thick] (a) to (b) to (15,0);
\draw [thick] (a) to (c) to (4.5,14.25);
\draw [thick] (a) to (d) to (-12,8.85);
\draw [thick] (a) to (e) to (-12,-8.85);
\draw [thick] (a) to (f) to (4.5,-14.25);

\draw [thick, dashed] (15,0) to (20,0);
\draw [thick, dashed] (4.5,14.25) to (6,19);
\draw [thick, dashed] (-12,8.85) to (-16,11.8);
\draw [thick, dashed] (-12,-8.85) to (-16,-11.8);
\draw [thick, dashed] (4.5,-14.25) to (6,-19);

\node at (-18,0) {${\mathcal B}_5$};


\node (b) at (50,0)  [shape=circle, fill=black]  {};
\node (c) at (43,9.5)  [shape=circle, fill=black] {};
\node (d) at (32,5.9) [shape=circle, fill=black] {};
\node (e) at (32,-5.9) [shape=circle, fill=black] {};
\node (f) at (43,-9.5) [shape=circle, fill=black] {};

\node (B) at (30,0)  [shape=circle, draw]  {};
\node (C) at (37,9.5)  [shape=circle, draw] {};
\node (D) at (48,5.9) [shape=circle, draw] {};
\node (E) at (48,-5.9) [shape=circle, draw] {};
\node (F) at (37,-9.5) [shape=circle, draw] {};

\draw [thick] (b) to (D) to (c) to (C) to (d) to (B) to (e) to (F) to (f) to (E) to (b);

\node at (58,0) {${\mathcal B}_5^{12}$};


\end{tikzpicture}

\end{center}
\caption{The dessins ${\mathcal B}_5$ (with a white vertex at infinity) and ${\mathcal B}_5^{12}$.}
\label{fig:B5}
\end{figure}

We will call $n$, $p$ and $d$ the {\sl valency}, the {\sl  field characteristic\/} and the {\sl dimension\/} of $\mathcal D$. (Of course, the field characteristic $p$ should not be confused with the Euler characteristic $\chi$ of $\mathcal D$.) Let $GP(n)$ and $GP(n,p)$ denote the set of all generalised Paley dessins which have valency $n$, and those which also have field characteristic $p$.
Thus $GP(1)$ consists of the primitive star dessins ${\mathcal St}_m$, those with $m=p$ prime. As another simple special case we have the following:

\medskip

\noindent{\bf Example 5} Let ${\mathcal D}\in GP(2)$, so $n=2$, $d=1$, $S=\{\pm 1\}\subseteq{\mathbb F}_p^*$ for some odd prime $q=p$, and $G$ is the dihedral group ${\rm D}_p$. If we choose $y\in T={\rm C}_p$ then $\mathcal D$ is the regular dessin ${\mathcal B}_p$ of type $(2,p,2)$ embedding a complete bipartite graph $K_{p,2}$ in $\Sigma$, with $p$ black vertices at the $p$th roots of unity and two white vertices at $0$ and $\infty$, each black and white pair joined by a single radial edge; the face-centres are obtained by multiplying the black vertices by $e^{\pi i/p}$. (This dessin is the union of ${\mathcal St}_p$ and its reflection in the equatorial plane of the sphere $\Sigma$.) If we choose $y\in G\setminus T={\rm D}_p\setminus{\rm C}_p$ then $\mathcal D$ is the dual dessin ${\mathcal B}_p^{12}$ of type $(2,2,p)$, with white vertices and face-centres transposed, and with edges around the unit circle. In either case the automorphism group $G$, generated by $z\mapsto e^{2\pi i/p}z$ and $z\mapsto\overline z$, acts primitively and faithfully on the black vertices.
Figure~\ref{fig:B5} shows planar drawings of ${\mathcal B}_5$, with a $5$-valent white vertex at infinity, and of ${\mathcal B}_5^{12}$.

\medskip

The next example is more instructive.

\medskip

\begin{figure}[h!]
\label{cube}
\begin{center}
\begin{tikzpicture}[scale=0.35, inner sep=0.7mm]

\node (a) at (0,0) [shape=circle, draw, fill=black] {}; 
\node (b) at (4,0) [shape=circle, draw, fill=black] {}; 
\node (d) at (2,3.46) [shape=circle, fill=black] {};  
\node (c) at (-2,3.46) [shape=circle, fill=black] {};  
\node (g) at (-4,0) [shape=circle, draw, fill=black] {}; 
\node (e) at (-2,-3.46) [shape=circle, fill=black] {};  
\node (f) at (2,-3.46) [shape=circle, fill=black] {};  

\node (u) at (6,3.46) [shape=circle, draw, fill=black] {}; 
\node (v) at (0,6.92) [shape=circle, draw, fill=black] {}; 
\node (w) at (-6,3.46) [shape=circle, fill=black] {};  
\node (x) at (-6,-3.46) [shape=circle, fill=black] {};  
\node (y) at (0,-6.92) [shape=circle, fill=black] {};  
\node (z) at (6,-3.46) [shape=circle, fill=black] {};  


\node at (0,-1) {$0$};
\node at (4,-1) {$1$};
\node at (2,2.46) {$4$};
\node at (-2,2.46) {$3$};
\node at (-4,-1) {$12$};
\node at (-2,-4.46) {$9$};
\node at (2,-4.46) {$10$};

\node at (6,2.46) {$5$};
\node at (0,5.92) {$7$};
\node at (-6,2.46) {$2$};
\node at (-6,-4.46) {$8$};
\node at (0,-7.92) {$6$};
\node at (6,-4.46) {$11$};


\draw [dashed] (2,8.58) to (-6,5.77 ) to (-8 , -2.3) to (-2,-8.58) to (6,-5.77 ) to (8,2.3) to (2,8.58);

\draw (7.5,0) to (-7.5,0); %
\draw ( 7,3.46) to (-6.65, 3.46); %
\draw ( -7,-3.46) to (6.65,-3.46); %
\draw ( -2.8,6.92) to (3.65,6.92); %
\draw ( 2.8,-6.92) to (-3.65,-6.92); %
\draw (3.85,6.75) to (-3.85,-6.75); %
\draw (3.8,-6.65) to (-3.8,6.65); %

\draw (-0.5,7.7) to (6.4,-4.1); %
\draw (0.5,-7.7) to (-6.4,4.1);  %
\draw (0.65,8.1) to (-6.4,-4.1);  %
\draw (-0.65,-8.1) to (6.4,4.1); %
\draw (4.7,5.7) to (7.6,0.8); %
\draw (-4.7,-5.7) to (-7.6,-0.8); %
\draw (-4.2,6.4) to (-7,1.8); 
\draw (4.2,-6.4) to (7,-1.8); 

 
  \end{tikzpicture}

\end{center}
\caption{A torus embedding of the Paley graph $P_{13}$.}
\label{fig:P13}
\end{figure}
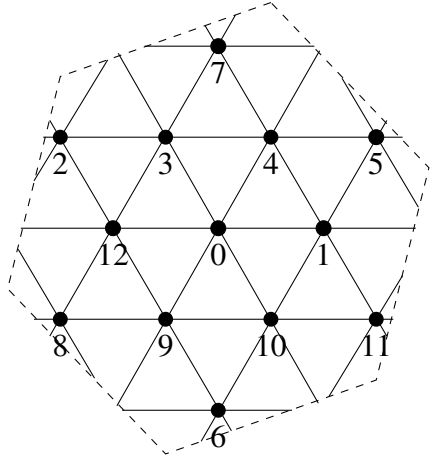

\noindent{\bf Example 6}  Let $q=13$ and $n=6$, so that $S$ is the group of quadratic residues (non-zero squares) in ${\mathbb F}_{13}$. This is generated by the element $4$, so we can define $x:t\mapsto 4t$. If we take $y: t\mapsto -t+b$ for any $b\ne 0$ in ${\mathbb F}_{13}$ then $x$ and $y$ generate the unique subgroup
$G={\rm AHL}_1(13)\cong{\rm C}_{13}\rtimes{\rm C}_6$
 of index $2$ in ${\rm AGL}_1(13)$. Since $y$ and $z:=(xy)^{-1}$ have orders $2$ and $3$ we obtain an orientably regular dessin $\mathcal D$ of type $(6,2,3)$ with automorphism group $G$. Since the white vertices all have valency $2$ we can ignore them and regard $\mathcal D$ as a $6$-valent triangular map $\mathcal M$. This map has $13$ vertices, $39$ edges and $26$ faces, so it has genus $1$; it is shown in Figure~\ref{fig:P13}, where opposite sides of the outer hexagon are identified to form a torus. The vertices of $\mathcal M$ are labelled with the elements of ${\mathbb F}_{13}$, and the neighbours of $0$ are successive powers of the generator $4$ of $S$. This map ($\{3,6\}_{3,1}$ in the notation of Coxeter and Moser~\cite[Section 8.4]{CM}) is chiral, with its mirror image formed similarly using the generator $4^{-1}=10$ of $S$; in each case, the embedded graph is the Paley graph $P_{13}$.

\medskip

The importance of generalised Paley dessins lies in Theorem~\ref{theo:main}, which asserts that
a regular dessin $\mathcal D$ is primitive and faithful if and only if it is a generalised Paley dessin.
Before proving this in Section~\ref{sec:proof}  we consider the basic properties of these dessins.

Let us define $\mathcal D$ to have {\sl colour constant\/} $c=c(x,y)=c({\mathcal D})\in{\mathbb Z}_n$ if $y\in Tx^c$, and for each $c\in{\mathbb Z}_n$ let $GP(n,p,c)$ denote the set of dessins ${\mathcal D}\in GP(n,p)$ with $c({\mathcal D)}=c$. If we take $c=0$, that is, if we choose $y$ to be a non-identity element of $T$ (as in ${\mathcal B}_p$, see Example~4), then $y$ has order $p$ and $z:=(xy)^{-1}$ has order $n$, so $\mathcal D$ has type $(n,p,n)$. In this case there are $q=p^d$ black vertices and faces, and $nq/p=np^{d-1}$ white vertices, so $\mathcal D$ has Euler characteristic
\[\chi=nq\left(\frac{2}{n}+\frac{1}{p}-1\right)=p^{d-1}(2p+n-np)\]
and genus $g=1-\frac{\chi}{2}$.
If we take $c=-1$ then $y$ has order $n$ while $z$, an element of $T$, has order $p$, so we obtain a dessin of type $(n,n,p)$, with $\chi$ and $g$ the same as for $c=0$ (see ${\mathcal B}_p^{12}$ in Example~4).
However, if we take some $c\ne 0, -1$ then $y$ has the same order as $x^c$, namely $m=n/\gcd(n,c)$, while $z$ has order $l=n/\gcd(n,c+1)$, so $\mathcal D$ has type $(n,m,l)$ with
\[\chi=nq\left(\frac{1}{n}+\frac{1}{m}+\frac{1}{l}-1\right)=q\left(1+\gcd(n,c)+\gcd(n,c+1)-n\right)\]
and $g=1-\frac{\chi}{2}$.

If $m=2$ then $\mathcal D$ can be regarded as a map $\mathcal M$  by ignoring its white vertices. This happens if and only if either $c=0$ and $p=2$, so that $\mathcal D$ has type $(n,2,n)$ and genus $1+2^{d-2}(n-4)$ with $n$ odd, or $c=n/2$ with $n$ even, so that $p$ is odd and $\mathcal D$ has type $(n,2,l)$ with $l=n$ or $n/2$ as $n\equiv 0$ or $2$ mod~$(4)$. In this last case the genus is $1+q(n-4)/4$ or $1+q(n-6)/4$ respectively.

Two generalised Paley dessins $\mathcal D$ and ${\mathcal D}'$ in $GP(n,p)$ are isomorphic if and only if their standard generating pairs $x, y$ and $x',y'$ are equivalent under ${\rm Aut}\,G$. This group is isomorphic to ${\rm A\Gamma L}_1(q)={\rm AGL}_1(q)\rtimes{\rm Gal}\,{\mathbb F}_q$, acting by conjugation on its normal subgroup $G$, regarded as a subgroup of ${\rm AGL}_1(q)$. It follows that $\mathcal D\cong{\mathcal D}'$ if and only if
\begin{itemize}
\item $x'=x^{\gamma}$ for some $\gamma\in{\rm Gal}\,{\mathbb F}_q\cong{\rm C}_d$;
\item $c(x',y')=c(x,y)$, i.e. $y'\in T(x')^c$ where $y\in Tx^c$.
\end{itemize}
If we denote $\mathcal D$ by ${\mathcal P}(x,c)$ then we have ${\mathcal P}(x,c)\cong{\mathcal P}(x',c)$ if and only if $x$ and $x'$ are Galois conjugate and $c=c'$.
It follows from condition~(b) of Lemma~\ref{lem:affineprim} that ${\rm Gal}\,{\mathbb F}_q$ acts semiregularly on generators $x$ for $S$, so it has $\phi(n)/d$ orbits on them; independently there are $n$ choices for $c$, so  we have:

\begin{prop}\label{prop:counting}
For each integer $n\ge 1$ and prime $p$ coprime to $n$, we have $|GP(n,p)|=n\phi(n)/d$, with $|GP(n,p,c)|=\phi(n)/d$ for each $c\in{\mathbb Z}_n$, where $d$ is the multiplicative order of $p$ {\rm mod}~$(n)$.
\end{prop}

The mirror image of ${\mathcal D}={\mathcal P}(x,c)$ is the generalised Paley dessin $\overline{\mathcal D}={\mathcal P}(x^{-1},c)$, in the same set $GP(n,p,c)$ as $\mathcal D$. Then $\mathcal D$ is real if $\overline{\mathcal D}\cong\mathcal D$, that is, the corresponding hypermap is fully regular; this happens if and only if $x$ is inverted by some $\gamma\in{\rm Gal}\,{\mathbb F}_q$, and since ${\rm Gal}\,{\mathbb F}_q$ is generated by the Frobenius automorphism $t\mapsto t^p$ we have:

\begin{prop}\label{prop:real}
The dessins ${\mathcal D}\in GP(n,p)$ are all real or all form chiral pairs as $-1$ is or is not congruent to a power of $p$ {\rm mod}~$(n)$.
\end{prop}

Note that if $p^i\equiv -1$ mod~$(n)$, then $p^{2i}\equiv 1$ mod~$(n)$, so either $n=2$ and $p$ is odd, or $n>2$, $d$ is even and $i\equiv d/2$ mod~$(d)$. Let ${\mathbb Z}_n^*$ denote the multiplicative group of units mod~$(n)$, and let ${\mathbb Z}^*_{n,p}$ denote the subgroup of ${\mathbb Z}_n^*$ of order $d$ generated by (the congruence class of) $p$. Then Proposition~\ref{prop:real} states that the dessins in $GP(n,p)$ are real or chiral as $-1\in{\mathbb Z}^*_{n,p}$ or not.

In order to characterise the pairs $n, p$ with $-1\in{\mathbb Z}^*_{n,p}$, let $n_r=r^{e_r}$ ($r$ prime, $e_r\ge 1$) be the prime powers appearing in the factorisation of $n$, and let $d_r$ be the multiplicative order of $p$ mod~$(n_r)$. Then we have the following:

\begin{prop}\label{prop:congruence}
There is some $i$ such that $p^i\equiv-1$ {\rm mod}~$(n)$ if and only if the orders $d_r$ for primes $r>2$ dividing $n$ all have the same $2$-part $2^k\ge2$, where $p\equiv-1$ {\rm mod}~$(n_2)$ if $n_2\ge 2$, and $k=1$ if $n_2\ge4$. In these circumstances the congruence $p^i\equiv-1$ {\rm mod}~$(n)$ is satisfied by $i={\rm lcm}\,\{d_r/2\mid r>2\}$ if some prime $r>2$ divides $n$, and otherwise by $i=1$.
\end{prop}

For instance, if $p=2$ then these conditions are satisfied for $n=3$ and $n=5$, but not for $n=15$: the orders $d_3=2$ and $d_5=4$ are both even but have different $2$-parts. More explicitly, $2^i\equiv-1$ mod~$(3)$ if and only if $i$ is odd, whereas $2^i\equiv-1$ mod~$(5)$ if and only if $i\equiv 2$ mod~$(4)$, and these two conditions are incompatible. However, if $p=2$ and $n=99$ then $d_3=6$ and $d_{11}=10$ are both even and have the same $2$-part, so the conditions are satisfied with $i=15$.

The primality of $p$ is irrelevant here: as for many other properties of generalised Paley dessins, the conditions depend only on the congruence class of $p$ mod~$(n)$. Of course, Dirichlet's Theorem shows that every congruence class in ${\mathbb Z}_n^*$ contains infinitely many primes.

\medskip

To avoid a longer digression, the proof of Proposition~\ref{prop:congruence} is postponed to Section~\ref{sec:appendix}.

\medskip

\noindent{\bf Warning.} If any dessin is real then it admits a reflection, so the underlying Riemann surface is symmetric. However, the converse is not true: for example the dessins in Example~6 (one of them shown in Figure~\ref{fig:P13}) are chiral, whereas their Riemann surface, the hexagonal torus, is symmetric.


\section{Operations on dessins}\label{sec:operations}

For each permutation $\pi$ of $\{0, 1, 2\}$ of order $2$ or $3$ there is a duality or triality operation $D^{\pi}$ on dessins, permuting the sets $V_0$, $V_1$ and $V_2$ of black and white vertices and face-centres by renaming them. These form a group $\Omega\cong{\rm S}_3$ of operations $D^{\pi}:{\mathcal D}\mapsto{\mathcal D}^{\pi}$ on dessins, preserving their underlying surface and automorphism group, and permuting the three periods in their type. If the subgroup of $\Omega$ leaving $\mathcal D$ invariant has order $k$, then $\mathcal D$ lies in an orbit of $\Omega$ consisting of $6/k$ mutually non-isomorphic dessins. 

The duality operation $D^{12}: {\mathcal D}\mapsto{\mathcal D}^{12}$, transposing white vertices and faces, preserves each set $GP(n,p)$. Since ${\mathcal D}^{12}={\mathcal P}(x,c')$ where $1+c+c'=0$ in ${\mathbb Z}_n$, it follows that ${\mathcal D}\cong{\mathcal D}^{12}$ if and only if $n$ is odd and $c=(n-1)/2$.
The other two duality operations, transposing black vertices with white vertices or faces, do not preserve the set of all generalised Paley dessins. For instance, if ${\mathcal D}\in GL(n,p,c)$ then the black/white dual ${\mathcal D}^{01}$ of $\mathcal D$ is a generalised Paley dessin if and only if $c\in{\mathbb Z}_n^*$, in which case ${\mathcal D}^{01}\in GP(n,p,c^{-1})$ where $c^{-1}$ denotes the multiplicative inverse of $c$ mod~$(n)$. In particular, ${\mathcal D}\cong{\mathcal D}^{01}$ if and only if some $\alpha\in{\rm Aut}\,G={\rm A\Gamma L}_1(q)$ transposes the standard generators $x$ and $y$ of $G$. There are now two cases to consider. If $c=1$ then $x$ and $y$ have the form $t\mapsto at$ and $t\mapsto at+b$ for some $a\ne 1$ and $b\ne 0$, in which case conjugation by $t\mapsto -t+b(1-a)^{-1}$ transposes them.
If $c\ne1$ then $x$ and $y$ are transposed if and only if $x$ and $x^c$ are transposed by some element of ${\rm Gal}\,{\mathbb F}_q$; this is equivalent to $d$ being even and $p^{d/2}\equiv c$ mod~$(n)$. Similarly, ${\mathcal D}\cong{\mathcal D}^{02}$ if and only if either $c'=1$ (equivalently $c=-2$), or $d$ is even and $p^{d/2}\equiv c'\;(=-1-c)$ mod~$(n)$.

The mutually inverse triality operations $D^{012}$ and $D^{210}$ act by cyclically permuting the standard generators $x, y, z$ of $G$, so $\mathcal D$ is invariant under one (and hence both) of these operations if and only if there is some $\alpha\in{\rm Aut}\,G$ inducing this $3$-cycle. If $c=1$ then the existence of $\alpha$ implies that $c'=1$, so that $n=3$; thus $d\le 2$, so $\alpha\in{\rm AGL}_1(q)$. Since $p$ is coprime to $n=3$, $\alpha$ must have the form $t\mapsto at+b$ where the element $a\in S$ has multiplicative order $3$. Conversely, given any $y\ne x$ in $Tx$ there exists such an automorphism $\alpha$ inducing the required $3$-cycle, so $\mathcal D$ is triality-invariant. If $c\ne 1$ then, arguing as for dualities, we find that $\mathcal D$ is triality-invariant if and only if $d$ is divisible by $3$ and $c\equiv p^{d/3}$ or $p^{2d/3}$ mod~$(n)$. To summarise, we have:

\begin{prop}\label{prop:Omega}
Let ${\mathcal D}\in GP(n,p,c)$. Then
\begin{itemize}
\item ${\mathcal D}\cong{\mathcal D}^{12}$ if and only if $n$ is odd and $c=(n-1)/2$;
\item ${\mathcal D}\cong{\mathcal D}^{01}$ if and only if either $c=1$, or $p^{d/2}\equiv c$ {\rm mod}~$(n)$ with $d$ even;
\item  ${\mathcal D}\cong{\mathcal D}^{02}$ if and only if either $c=-2$, or $p^{d/2}\equiv -1-c$ {\rm mod}~$(n)$ with $d$ even;
\item ${\mathcal D}\cong{\mathcal D}^{012}$, or equivalently ${\mathcal D}\cong{\mathcal D}^{210}$, if and only if either $n=3$ and $c=1$, or $p^{d/3}$ or $p^{2d/3}\equiv c$ {\rm mod}~$(n)$ with $d$ divisible by $3$.
\end{itemize}
\end{prop}

\begin{coro}\label{cor:Omega-inv}
A dessin ${\mathcal D}\in GP(n,p,c)$ is invariant under all six operations in $\Omega$ if and only if $n=3$ and $c=1$.
\end{coro}

The hole operations $H_j\;(j\in{\mathbb Z}_n^*)$ act on the set of all maps of valency $n$ by preserving the embedded graph and replacing the cyclic rotation of edges around each vertex (in the monodromy group) with its $j$th power. They preserve the automorphism group, but in general they may change the genus, since they change the faces. (For example, $H_2$ transposes the icosahedron and the great dodecahedron, which has pentagonal faces and has genus $4$.) One can apply these operations to regular dessins of type $(n,m,l)$, regarded as maps by ignoring their black and white vertex-colours, provided $j$ is coprime to both $n$ and $m$. This applies to dessins ${\mathcal D}\in GP(n,p,c)$ with $c\ne 0$, since in this case $m$ divides $n$; the result is another dessin $H_j({\mathcal D})\in GP(n,p,c)$. If $c=0$, so that $\mathcal D$ has type $(n,p,n)$, one can instead apply $D^{12}\circ H_j\circ D^{12}$, since ${\mathcal D}^{12}$ has type $(n,n,p)$. If $x$ and $x'$ are both generators of $S$ then $x'=x^j$ for some $j$ coprime to $n$, so we have:

\begin{prop}\label{prop:Hj}
The group of operations $H_j$ ($j\in{\mathbb Z}_n^*$) acts transitively on each set $GP(n,p,c)$.
\end{prop}

The stabiliser of each dessin, and hence the kernel of this action, is the subgroup $\{H_j\mid j\in{\mathbb Z}^*_{n,p}\}$ of order $d$ generated by $H_p$.

The most important hole operation, which applies for all valencies $n$, is $H_{-1}$, replacing each map or dessin with its mirror image by reflecting it in an edge. Since this operation has order $2$ and commutes with those in $\Omega$, we obtain a group $\Omega^*=\langle\Omega, H_{-1}\rangle\cong{\rm S}_3\times{\rm C}_2\cong{\rm D}_6$ of operations on maps and dessins.
If the subgroup of $\Omega^*$ leaving $\mathcal D$ invariant has order $k^*$, then $\mathcal D$ lies in an orbit of $\Omega^*$ consisting of $12/k^*$ mutually non-isomorphic dessins. In particular, if $\mathcal D$ is real (that is, $H_{-1}({\mathcal D})\cong{\mathcal D}$), then so are all other dessins in that particular orbit of $\Omega$, so $k^*=2k$ and the orbits of $\Omega$ and $\Omega^*$ are identical; otherwise, $k^*=k$ and the orbit of $\Omega^*$ is the disjoint union of that of $\Omega$ and its mirror image. Using Proposition~\ref{prop:real} and Corollary~\ref{cor:Omega-inv} we have:

\begin{coro}\label{cor:Omega*-inv}
A dessin ${\mathcal D}\in GP(n,p,c)$ is invariant under all twelve operations in $\Omega^*$ if and only if $n=3$, $p\equiv 2$ {\rm mod}~$(3)$ and $c=1$.
\end{coro}

\begin{figure}[h!]
\label{fig:cube}
\begin{center}
\begin{tikzpicture}[scale=0.05, inner sep=0.7mm]

\node (a3) at (-100,0) [shape=circle, draw] {};
\node (b3) at (-110,17.5)  [shape=circle, fill=black]  {};
\node (c3) at (-130,17.5)  [shape=circle, draw] {};
\node (d3) at (-140,0) [shape=circle, fill=black] {};
\node (e3) at (-130,-17.5) [shape=circle, draw] {};
\node (f3) at (-110,-17.5) [shape=circle, fill=black] {};

\node (A3) at (-80,0) [shape=circle, fill=black] {};
\node (B3) at (-100,35)  [shape=circle, draw]  {};
\node (C3) at (-140,35)  [shape=circle, fill=black] {};
\node (D3) at (-160,0) [shape=circle, draw] {};
\node (E3) at (-140,-35) [shape=circle, fill=black] {};
\node (F3) at (-100,-35) [shape=circle, draw] {};

\draw[thick] (a3) to (b3) to (c3) to (d3) to (e3) to (f3) to (a3);

\draw[thick] (a3) to (A3);
\draw[thick] (b3) to (B3);
\draw[thick] (c3) to (C3);
\draw[thick] (d3) to (D3);
\draw[thick] (e3) to (E3);
\draw[thick] (f3) to (F3);
\draw[thick, dashed] (A3) to (B3) to (C3) to (D3) to (E3) to (F3) to (A3);


\node (a1) at (20,0) [shape=circle, draw] {};
\node (b1) at (10,17.5)  [shape=circle, fill=black]  {};
\node (c1) at (-10,17.5)  [shape=circle, draw] {};
\node (d1) at (-20,0) [shape=circle, fill=black] {};
\node (e1) at (-10,-17.5) [shape=circle, draw] {};
\node (f1) at (10,-17.5) [shape=circle, fill=black] {};

\node (A1) at (40,0) [shape=circle, fill=black] {};
\node (B1) at (20,35)  [shape=circle, draw]  {};
\node (C1) at (-20,35)  [shape=circle, fill=black] {};
\node (D1) at (-40,0) [shape=circle, draw] {};
\node (E1) at (-20,-35) [shape=circle, fill=black] {};
\node (F1) at (20,-35) [shape=circle, draw] {};

\draw[thick] (a1) to (b1) to (c1) to (d1) to (e1) to (f1) to (a1);
\draw[thick] (a1) to (A1);
\draw[thick] (b1) to (B1);
\draw[thick] (c1) to (C1);
\draw[thick] (d1) to (D1);
\draw[thick] (e1) to (E1);
\draw[thick] (f1) to (F1);
\draw[thick, dashed] (A1) to (B1) to (C1) to (D1) to (E1) to (F1) to (A1);

\node (O) at (0,0) [shape=circle, fill=red] {};
\node (A2) at (30,17.5) [shape=circle, fill=red] {};
\node (B2) at (0,35) [shape=circle, fill=red] {};
\node (C2) at (-30,17.5) [shape=circle, fill=red] {};
\node (D2) at (-30,-17.5) [shape=circle, fill=red] {};
\node (E2) at (0,-35) [shape=circle, fill=red] {};
\node (F2) at (30,-17.5) [shape=circle, fill=red] {};

\draw [thick] (a1) to (O) to (d1);
\draw [thick] (b1) to (O) to (e1);
\draw [thick] (c1) to (O) to (f1);
\draw [thick] (a1) to (A2) to (b1) to (B2) to (c1) to (C2) to (d1) to (D2) to (e1) to (E2) to (f1) to (F2) to (a1);

\end{tikzpicture}

\end{center}
\caption{Two views of an $\Omega^*$-invariant torus dessin.}
\label{fig:invariant}
\end{figure}
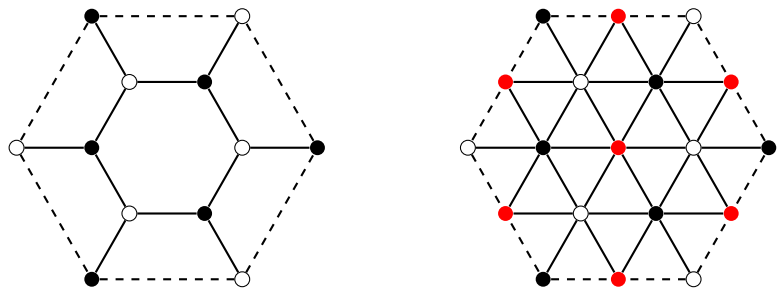

\noindent{\bf Example 7} For each prime $p\equiv 2$ mod~$(3)$ there is a unique dessin ${\mathcal D}\in GP(3,p,1)$. It has type $(3,3,3)$ and genus $1$. Since $d=2$ it has automorphism group $G\cong T\rtimes{\rm C}_3$ with $T\cong{\rm C}_p^2$. The smallest of these dessins, for $p=2$, is shown in Figure~\ref{fig:invariant}, on the left as a bicoloured map, and on the right, to emphasise its full symmetry under $\Omega^*$, as a tricoloured embedded graph, with red vertices at $\beta^{-1}(\infty)$; in each case, opposite sides of the outer hexagon are identified to form a torus. More generally, the triangle group $\Delta=\Delta(3,3,3)$, acting by isometries on $\mathbb C$, has a translation subgroup $\Delta_0\cong{\mathbb Z}^2$ of index $3$, and $\mathcal D$ is uniformised by the subgroup of index $p^2$ in $\Delta_0$ generated by its $p$th powers. For more on the dessins in $GP(3)$ see Example~9 in Section~\ref{sec:further}.

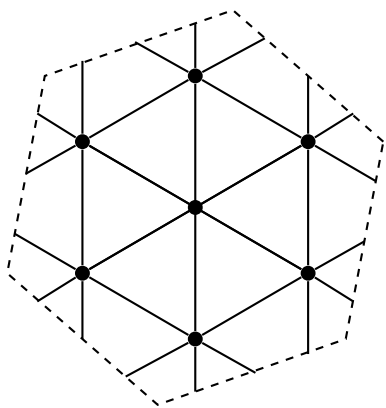
\begin{figure}[h!]
\label{fig:cube}
\begin{center}
\begin{tikzpicture}[scale=0.05, inner sep=0.7mm]

\draw[thick, dashed] (50,17.5) to (10,52.5) to (-40,35) to (-50,-17.5) to (-10,-52.5) to (40,-35) to (50,17.5); 

\node (A) at (30,17.5)  [shape=circle, fill=black]  {};
\node (B) at (0,35)  [shape=circle, fill=black]  {};
\node (C) at (-30,17.5)  [shape=circle, fill=black]  {};
\node (D) at (-30,-17.5)  [shape=circle, fill=black]  {};
\node (E) at (0,-35)  [shape=circle, fill=black]  {};
\node (F) at (30,-17.5)  [shape=circle, fill=black]  {};

\node (O) at (0,0) [shape=circle, fill=black] {};


\draw[thick] (A) to (D); 
\draw[thick] (0,49) to (B) to (E) to (0,-49);
\draw[thick] (30,35) to (A) to (F) to (30,-39);
\draw[thick] (-30,-35) to (D) to (C) to (-30,39);
\draw[thick] (C) to (F);
\draw[thick] (48,7) to (A) to (B) to (-16,44);
\draw[thick] (-48,-7) to (D) to (E) to (16,-44);
\draw[thick] (-45,9) to (C) to (B) to (18.5,45);
\draw[thick] (45,-9) to (F) to (E) to (-18.5,-45);
\draw[thick] (42,25) to (A) to (D) to (-42,-25);
\draw[thick] (-42,25) to (C) to (F) to (42,-25);

\end{tikzpicture}

\end{center}
\caption{A chiral dessin of type $(6,2,3)$, drawn as a map.}
\label{fig:623}
\end{figure}

Although these dessins ${\mathcal D}\in GP(3,p,1)$ for $p\equiv 2$ mod~$(3)$ are the only generalised Paley dessins which are invariant under the whole group $\Omega^*$, there are other dessins, such as the Fermat dessins, with this property. (Indeed, any characteristic subgroup of finite index in the free group $F_2$ uniformises such a dessin; for ${\mathcal F}_n$ this is the subgroup generated by the $n$th powers and the commutators.) At the other extreme, many generalised Paley dessins are invariant under only the identity operation in $\Omega^*$. These include the chiral torus dessins of type $(6,2,3)$ in $GP(6,p,3)$ for primes $p\equiv 1$ mod~$(6)$; see Figure~\ref{fig:623} for one of the chiral pair with $p=7$, drawn as a map.

One reason for the importance of understanding the action of this group $\Omega^*$ is that in Conder's census of regular and orientably regular maps and hypermaps~\cite{Conder}, each entry represents a single orbit of $\Omega^*$. In general, no information is given about the number of non-isomorphic objects in each orbit, although this can be deduced, after some computational labour, from the presentation for the automorphism group given there. The invariance results presented here provide an alternative solution to this problem for the generalised Paley dessins.

We have discussed invariance under operations in $\Omega$, and Proposition~\ref{prop:real} deals with $H_{-1}$, so there remain just the compositions of $H_{-1}$ with the duality and triality operations, of orders $2$ and $6$ respectively. For instance $H_{-1}({\mathcal D}^{12}) = (H_{-1}({\mathcal D}))^{12} = {\mathcal P}(x^{-1},c')$, and this is isomorphic to $\mathcal D$ if and only if  $-1\in{\mathbb Z}^*_{n,p}$ and $c=(n-1)/2$, with $n$ odd.

We have $H_{-1}({\mathcal D}^{01})\cong {\mathcal D}$ if and only if some $\alpha\in{\rm Aut}\,G$ transposes $x$ and $y^{-1}$. If $c=-1$, so that $y^{-1}\in Tx$, then conjugation by a suitable involution achieves this, as in the case of $D^{01}$ with $c=1$. If $c\ne -1$ it is necessary and sufficient that $p^{d/2}\equiv c^{-1}$ mod~$(n)$ with $d$ even. Similarly, $H_{-1}({\mathcal D}^{02})\cong {\mathcal D}$ if and only if either $c=0$, or $p^{d/2}\equiv (-1-c)^{-1}$ mod~$(n)$ with $d$ even. Finally, $\mathcal D$ is invariant under the mutually inverse operations $H_{-1}\circ D^{012}$ and $H_{-1}\circ D^{210}$ of order $6$ if and only if it is invariant under their cubes $H_{-1}$ and their squares $D^{012}$ and $D^{210}$, so by Propositions~\ref{prop:real} and \ref{prop:Omega} this happens if and only if (a) $-1\in{\mathbb Z}_{n,p}^*$, and either (b1) $n=3$ and $c=1$, or (b2) $p^{d/3}$ or $p^{2d/3}\equiv c$ {\rm mod}~$(n)$ with $d$ divisible by $3$. Now (a) and (b1) are equivalent to ${\mathcal D}\in GF(3,p,1)$ with $p\equiv 2$ mod~$(3)$, while (a) and (b2) are equivalent to $d$ being divisible by $6$ with $c\equiv p^{d/3}$ or $p^{2d/3}$ {\rm mod}~$(n)$ and $p^{d/2}\equiv -1$ mod~$(n)$. Of course, we have already met the first condition, that ${\mathcal D}\in GF(3,p,1)$ with $p\equiv 2$ mod~$(3)$ in Corollaries~\ref{cor:Omega-inv} and \ref{cor:Omega*-inv}, while the second condition is met when, for example, $n=7$, $p\equiv 3$ or $5$ mod~$(7)$ and $c=2$ or $4$.

This proves the following analogue of Proposition~\ref{prop:Omega}:

\begin{prop}\label{prop:Omega*}
Let ${\mathcal D}\in GP(n,p,c)$. Then
\begin{itemize}
\item ${\mathcal D}\cong H_{-1}({\mathcal D}^{12})$ if and only if $n$ is odd and $c=(n-1)/2$;
\item ${\mathcal D}\cong H_{-1}({\mathcal D}^{01})$ if and only if either $c=1$, or $p^{d/2}\equiv c^{-1}$ {\rm mod}~$(n)$ with $d$ even;
\item  ${\mathcal D}\cong H_{-1}({\mathcal D}^{02})$ if and only if either $c=0$, or $p^{d/2}\equiv (-1-c)^{-1}$ {\rm mod}~$(n)$ with $d$ even;
\item ${\mathcal D}\cong H_{-1}({\mathcal D}^{012})$, or equivalently ${\mathcal D}\cong H_{-1}({\mathcal D}^{210})$, if and only if either ${\mathcal D}\in GF(3,p,1)$ with $p\equiv 2$ {\rm mod}~$(3)$, or $d$ is divisible by $6$ with $c\equiv p^{d/3}$ or $p^{2d/3}$ {\rm mod}~$(n)$ and $p^{d/2}\equiv -1$ {\rm mod}~$(n)$.
\end{itemize}
\end{prop}

One can also ask which dessins ${\mathcal D}\in GP(n,p,c)$ are invariant under all the hole operations $H_j$ for $j\in{\mathbb Z}_n^*$. By Proposition~\ref{prop:Hj} this is equivalent to $GP(n,p,c)$ containing just one dessin (up to isomorphism), that is, to $\phi(n)/d=1$. Since $d$ is the multiplicative order of $p$ mod~$(n)$, this is equivalent to ${\mathbb Z}_n^*$ being a cyclic group with generator $p$. Now it is known (see~\cite[Theorem~6.11]{JJ}, for instance) that ${\mathbb Z}_n^*$ is cyclic if and only if $n=1$, $2$, $4$, $r^e$ or $2r^e$ where $r$ is an odd prime, in which case the number of generating elements is $\phi(\phi(n))$. Since we are assuming that $n\ge 2$ we have:

\begin{prop}\label{prop:Hj-inv}
The following are equivalent:
\begin{itemize}
\item A dessin ${\mathcal D}\in GP(n,p,c)$ is invariant under all the hole operations $H_j$ for $j\in{\mathbb Z}_n^*$;
\item $GP(n,p,c)$ contains only one dessin;
\item $p$ has multiplicative order $d=\phi(n$) {\rm mod}~$(n)$;
\item ${\mathbb Z}_n^*$ is a cyclic group with generator $p$;
\item $n=2$, $4$, $r^e$ or $2r^e$ where $r$ is an odd prime, and $p$ generates ${\mathbb Z}_n^*$.
\end{itemize}
\end{prop}

For instance, $n=18$ satisfies these conditions, with $\phi(\phi(18))=\phi(6)=2$ generating classes, namely $p\equiv 5$ and $p\equiv 11$ mod~$(18)$. Note that if these conditions are satisfied, they apply to all $n$ dessins in $GP(n,p)$.

Borrowing a term from the theory of regular maps, let us define a regular dessin to be {\sl kaleidoscopic\/} if it is invariant under $\Omega^*$ and all the relevant hole operations. From Corollary~\ref{cor:Omega*-inv} and Proposition~\ref{prop:Hj-inv} we have:

\begin{coro}
A dessin ${\mathcal D}\in GP(n,p,c)$ is kaleidoscopic if and only if $n=3$, $p\equiv 2$ {\rm mod}~$(3)$ and $c=1$.
\end{coro}

See Example~7 for more information on these dessins.


\section{Action on white vertices and face centres}\label{sec:white}

The automorphism group $G$ of any generalised Paley dessin $\mathcal D$ acts primitively on the black vertices. It acts primitively on the white vertices if and only if their stabilisers are maximal subgroups of $G$. If $c\ne 0$ then these stabilisers are the conjugates of $\langle x^c\rangle$, so this action is primitive if and only if $c$ is a unit mod~$(n)$. However, if $c=0$ then the stabilisers are cyclic subgroups of $T$, so the action is primitive if and only if $n$ is prime and $d=1$, that is, $p\equiv1$ mod~$(n)$.

The same argument applies to the face centres, except that we must replace $c$ with $c'=-1-c$, or equivalently with $c+1$. We therefore have:

\begin{prop}
Let $\mathcal D$ be a generalised Paley dessin, with automorphism group $G$. Then $G$ acts primitively on the white vertices (resp.~face-centres) of $\mathcal D$ if and only if either
\begin{itemize}
\item[(a)] $c$ (resp.~$c+1$) is a unit {\rm mod}~$(n)$, or
\item[(b)] $n$ is prime, $p\equiv 1$ {\rm mod}~$(n)$ and $c=0$ (resp.~$c=-1$).
\end{itemize}
\end{prop}

It follows that $G$ acts primitively on the white vertices and the face centres, in addition to the black vertices, if and only if either
\begin{itemize}
\item[(1)] $c$ and $c+1$ are both units {\rm mod}~$(n)$, or
\item[(2)] $n$ is prime, $p\equiv 1$ {\rm mod}~$(n)$ and $c=0$ or $-1$.
\end{itemize}
Condition~(1) implies that $n$ is odd, and conversely if $n$ is odd then one can satisfy this condition by taking $c=1$.


\section{Quotient dessins}\label{sec:quotients}

If $\mathcal D$ is a generalised Paley dessin with automorphism group $G$, then $\mathcal D$ has a quotient dessin ${\mathcal D}/T$ where $T$ is the translation subgroup of $G$, an elementary abelian group of order $q=p^d$. This dessin has a single black vertex since $T$ acts transitively on those of $\mathcal D$. Since $T$ is normal in $G$, ${\mathcal D}/T$ is a regular dessin, with automorphism group $G/T\cong S\cong{\rm C}_n$ where $n$ is the valency of $\mathcal D$. If $\mathcal D$ has colour constant $c\ne 0$ or $-1$ then the covering ${\mathcal D}\to{\mathcal D}/T$ is smooth, so ${\mathcal D}/T$ has the same type as $\mathcal D$ and has 
genus
\[\frac{1}{2}\left(1-\gcd(n,c)-\gcd(n,c+1)+n\right).\]
However, if $c=0$ or $-1$, so that $\mathcal D$ has type $(n,p,n)$ or $(n,n,p)$, then the covering is branched over the white vertices or face-centres; in these cases ${\mathcal D}/T$ has type $(n,1,n)$ or $(n,n,1)$ and is isomorphic to ${\mathcal St}_n^{01}$ or ${\mathcal St}_n^{02}$, each of genus $0$.

\medskip

\noindent{\bf Example 8} If $n=6$ then  ${\mathcal D}/T$ has genus $0, 2, 1, 1, 2, 0$ as $c=0, 1, \ldots, 5$ in ${\mathbb Z}_6$.

\medskip

In each case ${\mathcal D}/T$, as a regular dessin with an abelian automorphism group of exponent $n$, is a regular quotient of the Fermat dessin ${\mathcal F}_n$. Now ${\rm Aut}\,{\mathcal F}_n$ is isomorphic to ${\rm C}_n\times{\rm C}_n$, with standard generators $x$, $y$ and $z$ which act on points $[u,v,w]$ of the Fermat curve $u^n+v^n+w^n=0$ by multiplying $u$, $v$ or $w$ by $\zeta_n:=\exp(2\pi i/n)$ (so that $xyz=1$ since $u, v$ and $w$ are homogeneous coordinates). Then ${\mathcal D}/T$ is the quotient of ${\mathcal F}_n$ by the subgroup $\langle x^cy^{-1}, x^{c'}z^{-1}\rangle$ where $1+c+c'=0$ in ${\mathbb Z}_n$.


\section{Further examples}\label{sec:further}

We dealt with the case $n=2$ rather informally in Example~5. We can now use the above general arguments to deal with the next case.

\medskip

\noindent{\bf Example 9} Let $n=3$. If $p\equiv 1$ mod~$(3)$ then $d=1$, so there are $n\phi(n)/d=6$ dessins in $GP(3,p)$, forming one chiral pair for each $c\in{\mathbb Z}_3$, whereas if $p\equiv 2$ mod~$(3)$ then $d=2$, so there are $3$ real dessins, one for each $c\in{\mathbb Z}_3$. For $c=0$ or $2$ the dessins have type $(3,p,3)$ or $(3,3,p)$ respectively, and have genus $(p-1)/2$ or $(p-1)(p-2)/2$ as $p\equiv 1$ or $2$ mod~$(3)$, whereas for $c=1$ they have genus $1$ and type $(3,3,3)$. These torus dessins ${\mathcal T}_p$ of type $(3,3,3)$ are quotients of the universal hypermap of type $(3,3,3)$, obtained by factoring out normal subgroups $M$ of the triangle group $\Delta(3,3,3)$ which are contained in its translation subgroup with quotient ${\rm C}_p^d$; if $p\equiv 1$ mod~$(3)$, so that $d=1$, there are two such subgroups, conjugate in the extended triangle group, whereas if $p\equiv 2$ mod~$(3)$, so that $d=2$, there is one such subgroup, normal in it.  The self-dualities and -trialities are:
\begin{itemize}
\item ${\mathcal D}\cong{\mathcal D}^{12}$ if and only if $c=1$;
\item ${\mathcal D}\cong{\mathcal D}^{01}$ if and only if either $c=1$ or $p\equiv 2$ {\rm mod}~$(3)$ with $c=2$;
\item  ${\mathcal D}\cong{\mathcal D}^{02}$ if and only if either $c=1$ or $p\equiv 2$ {\rm mod}~$(3)$ with $c=0$;
\item ${\mathcal D}\cong{\mathcal D}^{012}\cong{\mathcal D}^{210}$ if and only if $c=1$.
\end{itemize}

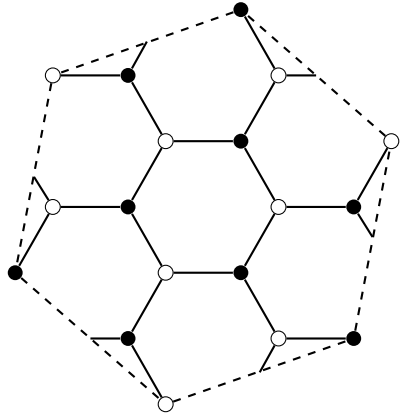
\begin{figure}[h!]
\label{fig:cube}
\begin{center}
\begin{tikzpicture}[scale=0.05, inner sep=0.7mm]

\node (a) at (20,0) [shape=circle, draw] {};
\node (b) at (10,17.5)  [shape=circle, fill=black]  {};
\node (c) at (-10,17.5)  [shape=circle, draw] {};
\node (d) at (-20,0) [shape=circle, fill=black] {};
\node (e) at (-10,-17.5) [shape=circle, draw] {};
\node (f) at (10,-17.5) [shape=circle, fill=black] {};

\node (A) at (40,0) [shape=circle, fill=black] {};
\node (B) at (20,35)  [shape=circle, draw]  {};
\node (C) at (-20,35)  [shape=circle, fill=black] {};
\node (D) at (-40,0) [shape=circle, draw] {};
\node (E) at (-20,-35) [shape=circle, fill=black] {};
\node (F) at (20,-35) [shape=circle, draw] {};

\node (a') at (50,17.5) [shape=circle, draw] {};
\node (b') at (10,52.5)  [shape=circle, fill=black]  {};
\node (c') at (-40,35)  [shape=circle, draw] {};
\node (d') at (-50,-17.5) [shape=circle, fill=black] {};
\node (e') at (-10,-52.5) [shape=circle, draw] {};
\node (f') at (40,-35) [shape=circle, fill=black] {};

\draw[thick] (a) to (b) to (c) to (d) to (e) to (f) to (a);
\draw[thick] (a) to (A) to (a');
\draw[thick] (b) to (B) to (b');
\draw[thick] (c) to (C) to (c');
\draw[thick] (d) to (D) to (d');
\draw[thick] (e) to (E) to (e');
\draw[thick] (f) to (F) to (f');
\draw[thick, dashed] (a') to (b') to (c') to (d') to (e') to (f') to (a'); 

\draw[thick] (A) to (45,-8);
\draw[thick] (B) to (30,35);
\draw[thick] (C) to (-15,44);
\draw[thick] (D) to (-45,8);
\draw[thick] (E) to (-30,-35);
\draw[thick] (F) to (15,-44);


\end{tikzpicture}

\end{center}
\caption{The dessin ${\mathcal T}_7$.}
\label{fig:T7}
\end{figure}

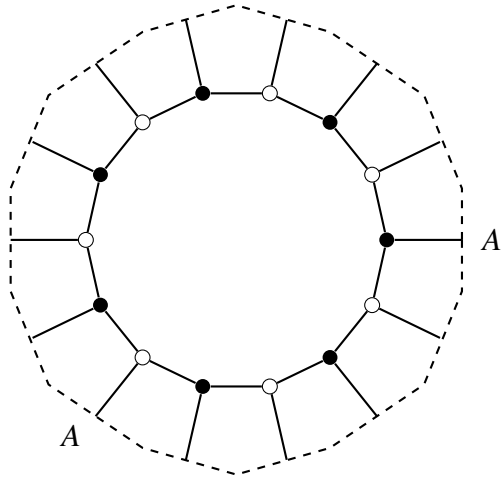
\begin{figure}[h!]
\label{fig:cube}
\begin{center}
\begin{tikzpicture}[scale=0.2, inner sep=0.7mm]


\node (a) at (10, 0) [shape=circle, fill=black] {};
\node (b) at (9.03, 4.34)  [shape=circle, draw]  {};
\node (c) at (6.23, 7.82)  [shape=circle, fill=black] {};
\node (d) at (2.23, 9.75) [shape=circle, draw] {};
\node (e) at (-2.23, 9.75) [shape=circle, fill=black] {};
\node (f) at (-6.23, 7.82)  [shape=circle, draw] {};
\node (g) at (-9.03, 4.34)  [shape=circle, fill=black]  {};
\node (h) at (-10, 0) [shape=circle, draw] {};
\node (i) at (-9.03, -4.34)  [shape=circle, fill=black]  {};
\node (j) at (-6.23, -7.82)  [shape=circle, draw] {};
\node (k) at (-2.23, -9.75) [shape=circle, fill=black] {};
\node (l) at (2.23, -9.75) [shape=circle, draw] {};
\node (m) at (6.23, -7.82)  [shape=circle, fill=black] {};
\node (n) at (9.03, -4.34)  [shape=circle, draw]  {};

\draw[thick] (a) to (b) to (c) to (d) to (e) to (f) to (g) to (h) to (i) to (j) to (k) to (l) to (m) to (n) to (a);

\draw[thick] (a) to (15,0);
\draw[thick] (b) to (13.54,6.52);
\draw[thick] (c) to (9.35,11.73);
\draw[thick] (d) to (3.35,14.63);
\draw[thick] (e) to (-3.35,14.63);
\draw[thick] (f) to (-9.35,11.73);
\draw[thick] (g) to (-13.54,6.52);
\draw[thick] (h) to (-15,0);
\draw[thick] (i) to (-13.54,-6.52);
\draw[thick] (j) to (-9.35,-11.73);
\draw[thick] (k) to (-3.35,-14.63);
\draw[thick] (l) to (3.35,-14.63);
\draw[thick] (m) to (9.35,-11.73);
\draw[thick] (n) to (13.54,-6.52);

\draw[thick, dashed] (15,0) to (15,3.42) to (12.48,9.62) to (6.22,13.84) to (0,15.61) to (-6.22,13.84) to (-12.48,9.62) to (-15,3.42) to (-15,0);
\draw[thick, dashed] (15,0) to (15,-3.42) to (12.48,-9.62) to (6.22,-13.84) to (0,-15.61) to (-6.22,-13.84) to (-12.48,-9.62) to (-15,-3.42) to (-15,0);

\node at (17,0) {$A$};
\node at (-11,-13) {$A$};

\end{tikzpicture}

\end{center}
\caption{A generalised Paley dessin of type $(3,3,7)$ and genus $3$.}
\label{fig:337}
\end{figure}

For instance, if $p=2$ we have the tetrahedron (with white vertices at the midpoints of its edges) for $c=0$ and its dual of type $(3,3,2)$ for $c=2$, together with a real torus dessin ${\mathcal T}_2$ of type $(3,3,3)$ for $c=1$, already shown in Figure~\ref{fig:invariant}. If $p=5$ we have real dessins of type $(3,5,3)$ and $(3,3,5)$ and genus $6$ (see RPH6.1 in~\cite{Conder} for the corresponding hypermaps) and a real torus dessin ${\mathcal T}_5$. However, if $p=7$ we have three chiral pairs, two of genus $3$ and types $(3,7,3)$ and $(3,3,7)$ on Klein's quartic curve (see CH3.1), together with a pair ${\mathcal T}_7$ and  $\overline{{\mathcal T}_7}$ of genus $1$; those dessins  of type $(3,3,7)$ and $(3,3,3)$ represent embeddings of the Heawood graph, with black and white vertices corresponding to the points and lines of the Fano plane ${\mathbb P}^2({\mathbb F}_2)$ (see~\cite{Sing}). The dessin ${\mathcal T}_7$ is shown in Figure~\ref{fig:T7}; with opposite sides of the outer hexagon  identified to form a torus. One of the chiral pair of type $(3,3,7)$ is shown in Figure~\ref{fig:337}; sides of the outer $14$-gon are identified, $A$ with $A$, and the others by ${\rm C}_7$ symmetry.

\medskip

The case $n=4$ is similar, so we omit it. Higher values of $n$ show a more interesting variety of phenomena.

\medskip

\noindent{\bf Example 10} Let $n=7$, so that $|GP(7,p)|$ is as follows:
\begin{itemize}
\item if $p\equiv 1$ mod~$(7)$ then $d=1$, so there are $42$ dessins;
\item if $p\equiv 6$ mod~$(7)$ then $d=2$, so there are $21$ dessins;
\item if $p\equiv 2$ or $4$ mod~$(7)$ then $d=3$, so there are $14$ dessins in each case;
\item  if $p\equiv 3$ or $5$ mod~$(7)$ then $d=6$, so there are $7$ dessins in each case.
\end{itemize}

If $p\equiv 3, 5$ or $6$ mod~$(7)$ then $p^{d/2}\equiv -1$ mod~$(7)$, so by Proposition~\ref{prop:real} the dessins ${\mathcal D}\in GP(7,p)$ are all real, whereas if $p\equiv 1, 2$ or $4$ mod~$(7)$ they form chiral pairs.
If $c=0$ or $6$ then $\mathcal D$ has type $(7,p,7)$ or $(7,7,p)$ respectively; in either case $\chi=p^{d-1}(7-5p)$, so $\mathcal D$ has genus
\[g=\frac{1}{2}\left(5p^d-7p^{d-1}+2\right).\]
However, if $c\ne 0$ or $6$ then $\mathcal D$ has type $(7,7,7)$ and genus
\[g=2p^d+1.\] 
Applying Proposition~\ref{prop:Omega} we have:
\begin{itemize}
\item ${\mathcal D}\cong{\mathcal D}^{12}$ if and only if $c=3$;
\item ${\mathcal D}\cong{\mathcal D}^{01}$ if and only if either $c=1$, or $c=6$ with $p\equiv 3, 5$ or $6$ mod~$(7)$; 
\item  ${\mathcal D}\cong{\mathcal D}^{02}$ if and only if either $c=5$, or $c=0$ with $p\equiv 3, 5$ or $6$ mod~$(7)$. 
\item ${\mathcal D}\cong{\mathcal D}^{012}\cong{\mathcal D}^{210}$ if and only if $c=2$ or $4$ with $p\equiv 2$, $3$, $4$ or $5$ mod~$(7)$.
\end{itemize}

Let us consider the six congruence classes $p\not\equiv 0$ mod~$(7)$ in turn. First let $p\equiv 1$ mod~$(7)$, so the $42$ dessins ${\mathcal D}\in GP(7,p)$ form $21$ chiral pairs, three for each $c\in{\mathbb Z}_7$. If $c=0$ or $6$ the dessins have type $(7,p,7)$ or $(7,7,p)$ and genus $g=5(p-1)/2$; these $12$ dessins, together with six duals of type $(p,7,7)$ which are not generalised Paley maps, form three orbits of length $6$ under the action of $\Omega^*$, with each dessin isomorphic to the mirror image its appropriate dual. For instance, if $p=29$ then $g=70$, and these three orbits correspond to items CH70.8, CH70.9 and CH70.10 in Conder's census~\cite{Conder}, each item representing one chiral pair of each of the above three types. (For details of the calculations supporting these and other similar claims, see the next section.) The remaining $30$ dessins in $GP(7,p)$ with $c=1, 2, 3, 4$ or $5$ all have type $(7,7,7)$ and genus $2p+1$; they form three orbits of length $6$, each orbit consisting of one chiral pair with each of $c=1,3$ and $5$, and there is one orbit of length $12$ consisting of three chiral pairs with each of $c=2$ and $4$. For instance, if $p=29$ these four orbits correspond to entries CH59.11--14 in~\cite{Conder}, the last having length $12$.

If $p\equiv 2$ or $4$ mod~$(7)$ then $d=3$, and in each case there are $14$ dessins, forming $7$ chiral pairs, one for each $c\in{\mathbb Z}_7$. If $c=0$ or $6$ the corresponding pair have type $(7,p,7)$ or $(7,7,p)$ and genus $g=(5p^3-7p^2+2)/2$; together with a chiral pair of type $(p,7,7)$ not in $GP(7,p)$ they form an orbit of length $6$. For instance, if $p=2$, so that $g=7$, this orbit corresponds to C7.2 in~\cite{Conder}, with the pair for $c=0$, regarded as maps, giving the Edmonds pair of orientably regular embeddings of the complete graph $K_8$. (Historically, this was the first discovery of a chiral pair of maps of genus $g\ne 1$.) If $c\ne 0$ or $6$ the corresponding pair have type $(7,7,7)$ and genus $g=2p^3+1$; these ten dessins form two orbits, of lengths $6$ (one pair each with $c=1, 3$ and $5$) and $4$ (with $c=2$ and $4$). For instance, if $p=2$ these orbits respectively correspond to CH17.7 and CH17.8 in~\cite{Conder}.

If $p\equiv 3$ or $5$ mod~$(7)$ then $d=6$, and in each case there are seven dessins, all real, one for each $c\in{\mathbb Z}_7$. If $c=0$ or $6$ they have type $(7,p,7)$ or $(7,7,p)$ and genus $(5p^6-7p^5+2)/2$; if $c\ne 0$ or $6$ they have type $(7,7,7)$ and genus $2p^6+1$. For each $p$ these seven dessins lie in three orbits of $\Omega$: one of length $2$ with $c=2$ or $4$, one of length $3$ with $c=1, 3$ and $5$, and another of length $3$ with $c=0$ and $6$ and a third dessin of type $(p,7,7)$ which is not a generalised Paley dessin. In all cases, even for $p=3$, the genus is too large for the corresponding hypermaps to appear in~\cite{Conder}.

Finally, if $p\equiv 6$ mod~$(7)$ then $d=2$ and there are $21$ dessins, all real, three for each $c\in{\mathbb Z}_7$, corresponding to the three mutually inverse pairs $x^{\pm 1}$ of generators of ${\mathbb Z}_7^*$. If $c=0$ or $6$ they have type $(7,p,7)$ or $(7,7,p)$ and genus $(5p^2-7p+2)/2$; if $c\ne 0$ or $6$ they have type $(7,7,7)$ and genus $2p^2+1$. There is one orbit of $\Omega$ of length $6$, consisting of the dessins with $c=2$ or $4$, and the other $15$ dessins, together with three more of type $(p,7,7)$, form six orbits of length $3$. Again the genus is too large for the corresponding hypermaps to be listed in~\cite{Conder}. 

\medskip

For comparison, we briefly consider an example where $n$ is composite, allowing a wider variety of types to appear:

\medskip

\noindent{\bf Example 11} Let $n=8$, so that $|GP(8,p)|$ is as follows:
\begin{itemize}
\item if $p\equiv 1$ mod~$(8)$ then $d=1$, with $32$ chiral dessins;
\item if $p\equiv 3$ or $5$ mod~$(8)$ then $d=2$, with $16$ chiral dessins in each case;
\item if $p\equiv 7$ mod~$(8)$ then $d=2$, with $16$ real dessins.
\end{itemize}

If $c=0$ or $7$ then $\mathcal D$ has type $(8,p,8)$ or $(8,8,p)$ respectively; in either case $\chi=2p^{d-1}(4-3p)$, so $\mathcal D$ has genus
\[g=1-\frac{\chi}{2}=3p^d-4p^{d-1}+1.\]
Thus if $p\equiv 1$ mod~$(8)$ then $d=1$ and $g=3(p-1)$; for instance, if $p=17$ then $g=48$ (see CH48.4--6 in~\cite{Conder}). If $p\equiv 3$, $5$ or $7$ mod~$(8)$ then $d=2$ and $g=(3p-1)(p-1)$; for instance, if $p=3$, $5$ or $7$ then $g=16$, $56$ or $120$  (see CH16.1 and CH56.5--6 in the first two cases). 

If $c=1$, $2$, $5$ or $6$ then $\mathcal D$ has type $(8,8,4)$ or $(8,4,8)$ as $c$ is odd or even. We have $\chi=-4p^d$, so $g=2p+1$ or $2p^2+1$ as $p\equiv 1$ mod~$(8)$ or otherwise. For example, if $p=17$ then $g=35$ (CH35.4--8); if $p=3$ then $g=19$ (CH19.8--9); if $p=5$ then $g=51$ (CH51.21-24); if $p=7$ then $g=99$ (RPH99.21--25). 

If $c=3$ or $4$ then $\mathcal D$ has type $(8,8,2)$ or $(8,2,8)$ respectively, with $\chi=-2p^d$ and hence $g=p+1$ or $p^2+1$ as $p\equiv 1$ mod~$(8)$ or otherwise. For example, if $p=17$ then $g=18$ (the chiral pair C18.1 of Paley maps); if $p=3$ then $g=10$ (the chiral pair C10.3 of self-dual Biggs embeddings of $K_9$); if $p=5$ then $g=26$ (C26.1); if $p=7$ then $g=50$ (R50.7).

\medskip

\noindent{\bf Example 12} The valencies $n=2, 3, 7$ and $8$ which we have considered in the above examples all have the property that $n+1$ is a prime power. More generally, if $n=p^e-1$ for some prime $p$, then for dessins ${\mathcal D}\in GP(n,p)$ we have $d=e$ and $G\cong{\rm AGL}_1(q)$ where $q=p^e$. Indeed, since ${\rm AGL}_1(q)$ acts doubly transitively on ${\mathbb F}_q$, but has no proper subgroups with this property, these are the only regular dessins in which $G$ is doubly transitive and faithful on the black vertices. If, in addition, either $c=0$ and $p=2$, or $n$ is even and $c=n/2$, then $m=2$ and $\mathcal D$ (with white vertices ignored) is a map $\mathcal M$, specifically one of the orientably regular embeddings of complete graphs $K_q$ constructed by Biggs in~\cite{Biggs}. For instance, if $n=10$ and $c=5$ we have the two chiral pairs C12.1 and C12.2 of Biggs embeddings $\mathcal M$ of $K_{11}$, corresponding to dessins $\mathcal D$ of type $(10,2,5)$ and genus $12$; if $n=15$ and $c=0$ we have the chiral pair C45.2 of Biggs embedding of $K_{16}$, corresponding to dessins of type $(15,2,15)$ and genus $45$. 

Similarly, if $2n=p^e-1$ for some prime $p$ (necessarily odd) then again dessins ${\mathcal D}\in GP(n,p)$ have $d=e$, but now $G$ is the subgroup ${\rm AHL}_1(q)$ of index $2$ in ${\rm AGL}_1(q)$. If $n$ is even, so that $q\equiv 1$ mod~$(4)$, and $c=n/2$, then $\mathcal D$ corresponds to one of the orientably regular Paley maps  $\mathcal M$ considered in~\cite{Jones}. For instance, if $n=6$ we obtain a chiral pair of torus embeddings of the Paley graph $P_{13}$, one of which was shown in Figure~\ref{fig:P13}.


\section{Typical calculations}\label{sec:calculations}

One can determine the values of $c$ for generalised Paley dessins appearing in~\cite{Conder} by abelianising the presentations of $G={\rm Aut}\,{\mathcal D}$ given there, thus determining the quotient dessins ${\mathcal D}/T$, which share these values of $c$. With $n$, $p$ and $c$ known, the invariance results in Section~\ref{sec:operations} then give the length of the orbit of $\Omega^*$ corresponding to each entry in~\cite{Conder}, that is, the number of non-isomorphic dessins it represents, information which is not given there.
 
\medskip

\noindent{\bf Example 13} As an illustration, the entries CH70.8, CH70.9 and CH70.10 in~\cite{Conder} represent the orbits of $\Omega^*$ containing orientably regular chiral hypermaps of genus $70$ and type $(7,7,29)$. In our terminology, these include the dessins $\mathcal D$ in $GP(7,29,c)$ with $c=0$ or $6$ (see Example~10 in Section~\ref{sec:further}).
In the case of CH70.8, the defining relations given for the automorphism group $G$ of $\mathcal D$, in terms of the standard generators $x$ and $y$ (or $R$ and $S$ in the notation of~\cite{Conder}), are:
 \[x^{-1}y^{-1}x^2y^2 = x^{-7} = y^{-1}x^{-1}yx^2yx^{-1}y^{-1} = 1.\]
 Abelianising these, we see from the first relation that $xy\in G'=T$, so $y\in Tx^{-1}$ and hence $c=-1=6$ in ${\mathbb Z}_7$. We obtain the same value for $c$ in the cases of CH70.9, where the relations are
 \[x^{-7} = xyx^3y^3 = y^{-1}x^{-1}yx^2yx^{-1}y^{-1} = y^{-7} = 1,\] 
and CH70.10, with relations
\[x^{-7} = y^{-1}x^{-1}yx^2yx^{-1}y^{-1} = x^{-1}y^2x^2y^{-1}xy = y^{-7} =1.\]
Since $c=6$ the invariance results in Section~\ref{sec:operations} show that the only operations in $\Omega^*$ preserving these dessins $\mathcal D$ of type $(7,7,29)$ are  $H_{-1}\circ{\mathcal D}^{01}$ and the identity, so the orbits CH70.8, CH70.9 and CH70.10 of $\Omega^*$ containing them have length $12/2=6$. Applying triality, we see that each of these orbits contains one chiral pair of generalised Paley dessins of each of the types $(7,7,29)$ and $(7,29,7)$ (the latter with $c=0$), together with one pair of type $(29,7,7)$ which are not generalised Paley dessins.
 
In order to determine which generalised Paley dessins $\mathcal D$ are in which of these orbits we need to know in each case how $x$ acts on $V$, or equivalently on $T$ by conjugation.  If we define $t:=xy$, so that $T=\langle t\rangle$, $y=x^{-1}t$ and $y^{-1}=t^{-1}x$, the second relation for CH70.8 becomes
 \[1 = x^{-1}.t^{-1}x.x^2.x^{-1}t.x^{-1}t = x^{-1}t^{-1}x^2tx^{-1}t = (t^{-1})^xt^{x^{-1}}t = t^{-r+r^{-1}+1},\]
 where we define $r\in{\mathbb F}_{29}$ by $t^x=t^r$. Thus $-r+r^{-1}+1=0$ in ${\mathbb F}_{29}$, so solving the equivalent quadratic equation $r^2-r-1=0$ gives $r=-5$ or $6$, and of these only $r=-5$ satisfies $r^7=1$, as implied by the second relation. (The third relation gives no useful information.) Thus the dessin $\mathcal D$ is defined within $GP(7,29,6)$ by the action $x: t\mapsto t^{-5}$ of $x$ on $T$. Since the multiplicative inverse of $-5$ in ${\mathbb F}_{29}$ is $-6$, in the mirror image $\overline{\mathcal D}=H_{-1}({\mathcal D})$ the action is given by $t\mapsto t^{-6}$. Applying the duality operation $D^{12}$ gives the second chiral pair of generalised Paley dessins in CH70.8, of type $(7,29,7)$ with $c=0$. Similar calculations show that in CH70.9 the corresponding chiral pair $\mathcal D$ and $\overline{\mathcal D}$ have $r=-4$ and $7$, and in CH70.10 they have $r=-9$ and $-13$. This accounts for all the elements $r$ of order $7$ in ${\mathbb F}_{29}^*$, and hence for all the $\phi(7)/d=6$ dessins in each of $GP(7,29,6)$ and $GP(7,29,0)$.
 
One can apply similar arguments to the dessins in $GP(7,29,c)$ for $c=1, 2, \ldots, 5$; these have type $(7,7,7)$ and genus $59$, and correspond to entries CH59.11--14 in~\cite{Conder}. (Again, see Example~10 in Section~\ref{sec:further}.) The dessin $\mathcal D$ representing CH59.11 has an automorphism group $G$ with defining relations
\[x^{-7} = x^{-1}y^{-1}x^{-2}yxy^{-1} = xy^{-1}x^{-2}y^{-3} =1.\]
Abelianising, we see that $y\in Tx^5$ and hence $c=5$. It follows that the only non-identity operation in $\Omega^*$ preserving $\mathcal D$ is $D^{02}$ so this orbit has length $12/2=6$. If we define $t:=yx^2\in T$, so that $y=tx^{-2}$ and $y^{-1}=x^2t^{-1}$, the second relation gives
\[1 = x^{-1}.x^2t^{-1}.x^{-2}.tx^{-2}.x.x^2t = xt^{-1}x^{-2}txt = (t^{-1})^{x^{-1}}t^xt = t^{-r^{-1}+r+1}\]
where $t^x=t^r$. Thus $-r^{-1}+r+1=0$ in ${\mathbb F}_{29}$, so $r^2+r-1=0$ giving $r=5$ or $-6$, and of these only $r=-6$ satisfies $r^7=1$. Thus in $\mathcal D$ the action of $x$ on $T$ is given by $t\mapsto t^{-6}$, and in $\overline{\mathcal D}$ by $t\mapsto t^{-5}$. The other four generalised Paley dessins in the orbit, with $c=1$ and $3$, are obtained by applying triality operations to this pair.

We obtain the same value for $c$, and hence for the length and structure of the orbit, in the cases of CH59.12, where the defining relations are
\[x^{-7} = xy^{-1}x^3yxy^{-1} = y^{-7} = x^{-1}yx^{-3}y^{-1}x^2y^{-1}= 1,\]
and CH59.13, where they are
\[x^{-7} = xyx^3y^{-2}x = yxyx^2y^{-1}xy = 1.\]
In these cases $x$ acts on $T$ in $\mathcal D$ and $\overline{\mathcal D}$ with $r=-4, 7$ and $-9, -13$ respectively.  (Note that we have now accounted for all six elements $r$ of order $n=7$ in ${\mathbb F}_{29}^*$, and hence for all nontrivial actions of $S$ on $T$.)

However, in CH59.14 the defining relations
 \[x^{-7} = xy^{-1}x^{-2}yx^{-1}y = y^2x^{-2}yxy = 1\]
imply that $c=2$. In this case no non-identity operations in $\Omega^*$ leave $\mathcal D$ invariant, so this orbit has length $12$, and all the duality and triality operations in $\Omega$ must be applied to the pair $\mathcal D$ and $\overline{\mathcal D}$ to find the other ten dessins, which have $c=1, 2$ or $4$. In $\mathcal D$ the action of $x$ on $T$ is given by $t\mapsto t^{-6}$, and in $\overline{\mathcal D}$ by $t\mapsto t^{-5}$.

In this example we have a favourable situation, with small parameters $n$, $p$ and especially $d$. In  general, if these are larger, and in particular if $T$ and $V$ have dimension $d>1$ over ${\mathbb F}_p$, then the calculations will be much more laborious, and the use of computer tools will be required.


\section{Proof of Theorem~\ref{theo:main}}\label{sec:proof}

In order to prove Theorem~\ref{theo:main} we need the concept of a {\sl Frobenius group}. This is a transitive permutation group in which stabiliser of any two points is the identity, but the stabiliser of any single point is not (that is, the group does not act regularly). Each Frobenius group $G$ is a semidirect product $N\rtimes G_0$ of a regular normal subgroup $N$ (the {\sl Frobenius kernel}, which is nilpotent by a theorem of Thompson~\cite{Thom}), and a point-stabiliser $G_0$ (a {\sl Frobenius complement}). For more background on Frobenius groups, see \cite[\S V.8]{Hup} or \cite{Pass}.

\begin{prop}\label{prop:primdes} Let $\mathcal D$ be a primitive regular dessin with automorphism group $G$ and black vertex set  $V$. Then
 the kernel $K$ of the action of $G$ on $V$ is a cyclic normal subgroup of $G$, and either
\begin{itemize}
\item[(a)] the permutation group $\overline{G}=G/K$ induced by $G$ on $V$ is regular of prime degree, or
\item[(b)] $\overline{G}$ s a Frobenius group, and $K$ is contained in the centre of $G$.
\end{itemize}
\end{prop}

\noindent{\sl Proof.} The stabilisers $G_v$ of vertices $v\in V$ are cyclic, and $K$ is contained in them, so $K$ is also cyclic. Since $\mathcal D$ is primitive, these stabilisers are maximal subgroups of $G$. Having the same stabiliser is an $G$-invariant equivalence relation, so either these stabilisers are all equal or they are mutually distinct. If they are all equal, and hence equal to $K$, then $\overline{G}$ is a regular permutation group on $V$, of prime degree by primitivity. We may therefore assume that distinct vertices $v, w\in V$ have distinct stabilisers. These stabilisers are abelian, so the centraliser $C_G(G_v\cap G_w)$ of their intersection contains both $G_v$ and $G_w$ and is therefore equal to $G$; thus $G_v\cap G_w$ is in the centre of $G$, so since the stabilisers $G_u$ ($u\in V$) are mutually conjugate $G_v\cap G_w$ fixes every black vertex and is therefore contained in the kernel $K=\cap_{u\in V}G_u$ of the action of $G$ on $V$. Clearly $G_v\cap G_w\ge K$, so $G_v\cap G_w = K$ and thus $K$ is contained in the centre of $G$. In the permutation group $\overline{G}=G/K$ induced by $G$ on $V$, the stabiliser $\overline{G}_v\cap \overline{G}_w=(G_v\cap G_w)/K$ of each pair $v\ne w$ in $V$ is the identity. Since $G_v>K$ we have $\overline{G}_v>1$, so $\overline{G}$ is a Frobenius group. \hfill$\square$

\medskip

\noindent{\bf Example 14}  To see how $K$ need not be a central subgroup in case~(a), let
\[G=\langle x,y\mid x^n=y^p=1,\, x^y=x^r\rangle \cong{\rm C}_n\rtimes{\rm C}_p,\]
where $r$ is an element of prime order $p$ in ${\mathbb Z}_n^*$. Then $G={\rm Aut}\,{\mathcal D}$ for a regular dessin $\mathcal D$, with standard generators $x$ and $y$. There are $p$ black vertices of valency $n$, and $n$ white vertices of valency $p$. Now $G$ permutes the black vertices primitively, inducing a regular permutation group $\overline{G}=G/K\cong{\rm C}_p$ on them, with a kernel $K=\langle x\rangle$ which is cyclic but not central in $G$. For instance, if $p=2$ and $r=-1$ then $G\cong{\rm D}_n$ and $\mathcal D$, regarded as a map, is the $n$-valent hosohedron (or beachball) on the sphere.

\medskip

\noindent{\sl Proof of Theorem~\ref{theo:main}.} By their construction, generalised Paley dessins are primitive and faithful, as are $p$-star dessins for each prime $p$.

Conversely, if a regular dessin $\mathcal D$ is primitive and faithful, then $G=\overline{G}$ and by Proposition~\ref{prop:primdes} this is a regular or Frobenius group. In the first case Lemma~\ref{lem:regonV} shows that $\mathcal D$ is a star map, so we may assume that $G$ is a Frobenius group.

As in any Frobenius group, the Frobenius kernel $N$ of $G$ is a nilpotent normal subgroup, acting regularly on $V$. Since $G$ acts primitively, $N$ is characteristically simple (for otherwise it would contain a non-identity normal and hence transitive subgroup of $G$ as a proper subgroup) and it is therefore a direct product is isomorphic simple groups. Since $N$ is nilpotent these must be cyclic groups of prime order, so $N$ is elementary abelian, of order $q=p^d$ for some prime~$p$, and we can regard $N$ as a $d$-dimensional vector space over the field ${\mathbb F}_p$. Since $N$ acts regularly on $V$ we can identify $V$ with $N$, acting on itself by translations. The Frobenius complement $G_0$, the stabiliser of the vector $0$, acts by conjugation on this vector space $V$ as a subgroup of ${\rm GL}_d(p)$, which is irreducible since $G$ is primitive. Since $G_0$ is abelian it follows from Schur's Lemma that one can identify $V$ with the additive group of the field ${\mathbb F}_q$, and $G_0$ with a subgroup $S$ of the multiplicative group ${\mathbb F}_q^*={\mathbb F}_q\setminus\{0\}$, which generates $V$ (see, for example, \cite[Satz II.3.10]{Hup} with $f=1$). Thus $G=V\rtimes G_0$ acts on $V$ as a primitive subgroup of ${\rm AGL}_1(q)$, as in Lemma~\ref{lem:affineprim}. The black vertices of $\mathcal D$ correspond to the cosets of the subgroup $G_0=S$, and the cyclic order of edges around the vertex $0$ is given by successive powers of a generator $x$ of $S$, so $\mathcal D$ is a generalised Paley dessin, as in Example~3. \hfill$\square$


\section{Galois action}\label{sec:Galois}

The absolute Galois group $\Gamma={\rm Gal}\,\overline{{\mathbb Q}/{\mathbb Q}}$ acts on isomorphism classes of dessins, preserving their automorphism group and type. It therefore acts on the set of generalised Paley dessins $\mathcal D$, preserving their valency $n$ and characteristic $p$. Recall that for a given $n$ and $p$ (not dividing $n$) there are $n\phi(n)/d$ dessins, with $\phi(n)/d$ for each colour constant $c\in{\mathbb Z}_n$, where $d$ is the multiplicative order of $p$ mod~$(n)$. In fact, $\Gamma$ also preserves the value of $c$. One way to see this is to note that there is an induced action of $\Gamma$ on the set of quotient dessins ${\mathcal D}/T$, since the translation group $T$ of order $q=p^d$ is a characteristic subgroup of $G:={\rm Aut}\,{\mathcal D}$. These are regular dessins with abelian automorphism groups (isomorphic to ${\rm C}_n$), so they are defined over $\mathbb Q$ (see~\cite{WM}) and hence invariant under $\Gamma$. Since the values of $c$ are in bijective correspondence with the isomorphism classes of these dessins, they too are preserved by $\Gamma$. Thus, for each triple $n$, $p$ and $c$ there is an action of $\Gamma$ on the set $GP(n,p,c)$ of generalised Paley dessins with these parameters.

\medskip

\noindent{\bf Question:} is there a more direct proof that $c$ is invariant under $\Gamma$?

\medskip

Since the group of hole operations $H_j$ ($j\in{\mathbb Z}_n^*$) acts transitively on each of these sets $GP(n,p,c)$, it follows from~\cite[Theorem~2]{JSW} that they are the orbits of $\Gamma$ on generalised Paley dessins, and that each dessin $\mathcal D$ in this orbit is defined over a subfield $K$ of the cyclotomic field ${\mathbb Q}(\zeta_n)$, where $\zeta_n:=e^{2\pi i/n}$: specifically, $K$ is the fixed field of the subgroup $H=\{j\in{\mathbb Z}_n^*\mid H_j({\mathcal D})\cong{\mathcal D}\}$, where we identify ${\rm Gal}\,{\mathbb Q}(\zeta_n)/{\mathbb Q}$ with ${\mathbb Z}_n^*$ in the usual way. (This generalises a similar result for generalised Paley maps proved in~\cite{Jones}; see also~\cite[Section~9.1.1]{JW}.) The length of this orbit is the degree $|K:{\mathbb Q}|=\phi(n)/d$ of this extension. Proposition~\ref{prop:Hj-inv} tells us when this degree is $1$, so we have:

\begin{prop}\label{prop:def/Q}
A dessin ${\mathcal D}\in GP(n,p,c)$ is defined over $\mathbb Q$ if and only if
$n=2$, $4$, $r^e$ or $2r^e$ where $r$ is an odd prime, and $p$ generates ${\mathbb Z}_n^*$.
\end{prop}


\section{Equations}\label{sec:equations}

A significant and difficult problem is that of converting combinatorial data about a dessin into defining equations. Even with the group-theoretic advantage of a regular dessin, this has been accomplished in only a few families of simple cases: see~\cite[Section 9.2]{JW}, for instance. One of these families consists of generalised Paley dessins.

\medskip

\noindent{\bf Example 15} Suppose that ${\mathcal D}\in GP(n,p)$ where $p\equiv 1$ mod~$(n)$, so that $d=1$ and $G$ is a metacyclic group ${\rm C}_p\rtimes{\rm C}_n$. If we take $c=0$ then $\mathcal D$ has type $(n,p,n)$ and genus $(n-2)(p-1)/2$; the embedded graph is the complete bipartite graph $K_{p,n}$, with $p$ black and $n$ white vertices of valencies $n$ and $p$. There are $\phi(n)$ such dessins $\mathcal D$ in $GP(n,p,0)$, one for each element $x$ of order $n$ in ${\mathbb F}_p^*$. As explained in Section~\ref{sec:Galois} these $\phi(n)$ dessins are permuted transitively by the hole operations $H_j$ ($j\in{\mathbb Z}_n^*$) and by the absolute Galois group, acting as ${\rm Gal}\,{\mathbb Q}(\zeta_n)/{\mathbb Q}\cong{\mathbb Z}_n^*$. It is therefore sufficient to give defining equations, with coefficients in ${\mathbb Q}(\zeta_n)$, for just one of these dessins, since equations for the others can then be obtained by Galois conjugacy. Here we are in the favourable situation that $\mathcal D$ is a cyclic covering of ${\mathcal D}/T$, which is a dessin on the Riemann sphere $\Sigma$. In fact, ${\mathcal D}/T={\mathcal St}_n^{01}$, with a black vertex of valency $n$ at $0$, and white vertices of valency $1$ at the $n$th roots of $1$, where the $p$-sheeted covering is branched. In Example~4 of~\cite{JSW} these facts are used to give (with a few minor changes of notation)  the affine model
\[w^p=\prod_{j=1}^n\left(z-\zeta_n^j\right)^{u^j},\]
in variables $w, z\in{\mathbb C}$, for one of the underlying curves, where $x$ acts on ${\mathbb F}_p$ by $t\mapsto t^x=ut$ for $u\in{\mathbb Z}_n^*$. The corresponding Bely\u\i\/ function $(w,z)\mapsto z^n$ is the composition of the covering ${\mathcal D}\to{\mathcal D}/T={\mathcal St}_n^{01}$, $(w,z)\mapsto z$, with the Bely\u\i\/ function $z\mapsto z^n$ of ${\mathcal St}_n^{01}$. Formulae for the standard generators of $G$ are given in~\cite{JSW}; the automorphisms of order $p$ are obvious, but those of order $n$ require careful calculation.
Replacing $\zeta_n$ with other primitive $n$th roots of $1$ gives the remaining curves and dessins in this Galois orbit. (This example is also considered in~\cite[Section~9.2.4]{JW}, where for simplicity $n$ is assumed to be prime.) One can extend this to the Galois orbit $GP(n,p,-1)$ by using the duality $D^{12}$, replacing the variable $z$ with $z/(z-1)$.

\medskip

\noindent{\bf Problem.} When $c=0$ or $-1$, so that ${\mathcal D}/T$ has genus $0$, can one extend these results for $d=1$ to cases where $d>1$? A simple example is the tetrahedron, the unique dessin ${\mathcal D}\in GP(3,2,0)$: this is on the sphere $\Sigma$, with $d=2$ and $G\cong{\rm A}_4\cong{\rm V}_4\rtimes{\rm C}_3$, and it has Bely\u\i\/ function
\[\beta(z)=-64\frac{(z^3-1)^3}{(z^3-8)^3z^3}.\]
(Magot and Zvonkin give this and many other Bely\u\i\/ functions on  the sphere in~\cite{MZ}.) Can one find defining equations for the unique dessin ${\mathcal D}\in GP(3,p,0)$ for other primes $p\equiv 2$ mod~$(3)$, which all give $d=2$ (and have genus $(p-1)(p-2)/2$ with $G\cong{\rm V}_{p^2}\rtimes{\rm C}_3$)?


\section{Non-faithful actions}\label{sec:non-faithful}

In Section~\ref{sec:examples} we saw some examples of non-faithful actions of $G$, on black vertices or more generally on all special points.  In the case of maps, Li and \v Sir\'a\v n~\cite{LS} have dealt thoroughly with non-faithful actions, and most of their results also apply to dessins, so here we just give a few examples to show how such actions on dessins are related to faithful actions.

If $\mathcal D$ is a regular dessin with automorphism group $G$ acting primitively but {\sl not\/} faithfully on its black vertex set $V$, then Proposition~\ref{prop:primdes} shows that the kernel $K$ of the action of $G$ on $V$, the intersection of the stabilisers $G_v$ $(v\in V$), is a cyclic normal subgroup of $G$, contained in the centre of $G$ if $G/K$ is a Frobenius group. Now $G/K$ acts primitively and faithfully on the vertex set $V/K=V$ of ${\mathcal D}/K$, so Theorem~\ref{theo:main} shows that ${\mathcal D}/K$ is a generalised Paley dessin.
Then $\mathcal D$ is a $k$-sheeted regular cyclic covering of ${\mathcal D}/K$, branched over its black vertices (and possibly its white vertices and/or face-centres), where $k=|K|$. In Example~2, for instance, where $\mathcal D$ is the Fermat dessin ${\mathcal F}_n$, the kernel $K$ is the subgroup $\langle x\rangle\cong{\rm C}_n$ of $G={\rm C}_n\times{\rm C}_n$ fixing the $n$ black vertices and cyclically permuting the $n$ white vertices and the $n$ faces, so that ${\mathcal D}/K$ is the star dessin ${\mathcal St}_n$ (see Example~1); the $n$-sheeted regular covering ${\mathcal D}\to{\mathcal D}/K$ is branched over the black vertices of ${\mathcal St}_n$, but not over its white vertex or face centre.
In Examples~3 and 4 the kernel $K$ is the commutator subgroup $G'$ of $G$, and ${\mathcal D}/K$ is a primitive Fermat dessin ${\mathcal F}_p$; in this case the $p$-sheeted covering is branched over all $3p$ special points of ${\mathcal F}_p$.

One can reverse this process, and construct non-faithful primitive dessins $\tilde{\mathcal D}$ from a faithful primitive dessin $\mathcal D$. Let ${\mathcal D}\in GP(n,p, c)$, with automorphism group $G=T\rtimes S$ as before. Let $\tilde S$ be a cyclic group of order $kn$ for some $k\in{\mathbb N}$, so that there is an epimorphism $\theta: \tilde{S}\to S$, and let ${\tilde G}=T\rtimes\tilde{S}$ where the action of $\tilde S$ on $T$ is given by composing that of $S$ with $\theta$. The unique subgroup $K=\ker\theta$ of order $k$ in $\tilde S$ is then a central cyclic subgroup of $\tilde G$, and $\theta$ extends to an epimorphism $\tilde{\theta}:\tilde{G}\to G$ with kernel $K$, so that $\tilde{G}/K\cong G$. If we choose a generator $\tilde x$ of $\tilde S$ with $\theta(\tilde{x})=x$, and any element $\tilde{y}\in\tilde{\theta}^{-1}(y)\setminus\tilde{S}$, they generate $\tilde G$ and hence determine a primitive regular dessin $\tilde{\mathcal D}$ which is a cyclic covering of $\mathcal D$, with $K$ as  the kernel of the action of $\tilde{G}={\rm Aut}\,\tilde{\mathcal D}$ on the black vertices. There are $q$ of these, of valency $kn$, while the numbers and valencies of the white vertices and faces, and hence the type and genus of $\tilde{\mathcal D}$, depend on the choice of $\tilde y$.

\medskip

\noindent{\bf Example 16} If $k$ is coprime to the order $nq$ of $G$ then the extension $\tilde G$ of $G$ by $K$ splits, by the Schur--Zassenhaus Theorem. Since $K$ is in the centre of $\tilde G$ we have $\tilde{G}=G\times K$, so that $\tilde{D}$ is the direct product (minimal common cover) of $\mathcal D$ and a regular dessin $\mathcal K$ with one black vertex and with automorphism group $K\cong{\rm C}_k$. There are $k$ such dessins $\mathcal K$, one for each choice of $y$ as a power of the standard generator $x$ of $K$; for instance, taking $y=1$ gives ${\mathcal K}={\mathcal St}_k^{01}$.

\medskip

\noindent{\bf Example 17} For any odd prime $p$ the unique dessin $\mathcal D$ in $GP(2,p,0)$ is the spherical dessin ${\mathcal B}_p$ of type $(2,p,2)$ in Example~5, with automorphism group $G\cong{\rm D}_p$ acting primitively and faithfully on the $p$ black vertices (and on the $p$ faces). Here $S\cong{\rm C}_2$, so we can take $\tilde S$ to be a cyclic group of any even order $2k>2$, with a unique epimorphism $\theta:\tilde{S}\to S$. To obtain a covering $\tilde{\mathcal D}\to{\mathcal D}$ we need to choose $\tilde{y}\in T\tilde{x}^{\tilde c}\setminus \tilde{S}$ for some even $\tilde{c}\in{\mathbb Z}_{2k}$. For instance, if $\tilde{c}=0$ then $\tilde{y}\in T$, so $\tilde{\mathcal D}$ has type $(2k,p,2k)$ and genus $p-1$; there is branching over the black vertices and face-centres, but not the white vertices. For $p=3, 5, 7, 11$ and $13$ the corresponding hypermap in~\cite{Conder} is RPH2.2, RPH4.9, RPH6.5, RPH 10.18 and RPH12.4. If $\tilde{c}\ne0$ then $\tilde{\mathcal D}$ has type $(2k,\tilde{m},\tilde{l})$ where $\tilde{m}$ and $\tilde{l}$ are the orders of $\tilde{y}$ and $\tilde{}z$;
for instance, if $\tilde{c}=2$, so that $\tilde{y}=t\tilde{x}^2=\tilde{x}^2t$ for some non-identity $t\in T$, then $\tilde{m}={\rm lcm}\,(k,p)$.

\section{Appendix: the proof of Proposition~\ref{prop:congruence}}\label{sec:appendix}

Clearly $p^i\equiv -1$ mod~$(n)$ if and only if $p^i\equiv -1$ mod~$(n_r)$ for each prime power $n_r=r^{e_r}$ ($r$ prime) appearing in the factorisation of $n$.

If $r>2$ then ${\mathbb Z}_{n_r}^*$ is cyclic, and $-1$ is its unique involution, so $p^i\equiv -1$ mod~$(n_r)$ if and only if $i\equiv d_r/2$ mod~$(d_r)$ where $p$ has even order $d_r$ in ${\mathbb Z}_{n_r}^*$. If $r=2$ then the same applies provided $n_2=4$, but if $n_2>4$ then ${\mathbb Z}_{n_2}^*=\langle 5\rangle\times\langle -1\rangle\cong{\rm C}_{n_2/4}\times{\rm C}_2$ (see~\cite[Theorem~6.10]{JJ}, for instance), a non-cyclic group with three involutions, $\pm 5^{n_2/8}=n_2/2\pm1$ and $-1$. In this case, and also when $n_2=4$, we have $p^i\equiv -1$ mod~$(n_2)$ if and only if $p\equiv -1$ mod~$(n_2)$ and $i$ is odd;  when $n_2=2$ one can omit the requirement that $i$ is odd. To summarise, $p^i\equiv -1$ mod~$(n)$ if and only if
\begin{itemize}
\item [a)] $i\equiv d_r/2$ mod~$(d_r)$ where $p$ has even order $d_r$ in ${\mathbb Z}_{n_r}^*$, for each odd $r$, and
\item [b)] $p\equiv -1$ mod~$(n_2)$ if $n$ is even, with any $i$ or any odd $i$ as $n_2=2$ or $n_2>2$.
\end{itemize}
The congruences in (a) require that $i$ is an odd multiple of $d_r/2$ for each odd $r$, so they are satisfied by the same $i$ if and only if all such $d_r$ have the same $2$-part $2^k\ge2$ in their factorisation: one can then take $i$ to be the least common multiple of the integers $d_r/2$ for odd $r$. If $n\equiv 0$ mod~$(4)$ then (b) implies that $i$ must be odd, so by (a) we have $d_r\equiv 2$ mod~$(4)$ for each odd $r$, that is, $k=1$.
If $n\not\equiv 0$ mod~$4$ then (b) imposes no restrictions on $p$, apart from being coprime to $n$.

\newpage


\centerline{\bf Acknowledgments}

\medskip

The authors are grateful to Rub\'en Hidalgo for some very helpful comments and corrections. The second author acknowledges support from the APVV Research Grant APVV-19-0308 and
from the VEGA Research Grants 1/0423/20 and 1/0727/22.


\end{document}